\documentclass{article}
\usepackage{amsthm}
\usepackage{amsfonts}
\usepackage{amssymb}
\usepackage{amsgen}
\usepackage{amsmath}
\usepackage{amsopn}
\usepackage{verbatim}
\usepackage{xypic}
\usepackage{xspace}
\usepackage{multicol}
\usepackage{makeidx}
\usepackage{eepic}
\usepackage{upref}

\theoremstyle{plain}
\newtheorem{thm}{Theorem}[section]
\newtheorem{lem}[thm]{Lemma}
\newtheorem{prop}[thm]{Proposition}
\newtheorem{cor}[thm]{Corollary}

\theoremstyle{definition}

\newtheorem{rem}[thm]{Remark}

\newtheorem{eg}[thm]{Example}

\newcommand{\ors}{{\vec s\,}}

\newcommand{\Z}{{\mathbb Z}}

\newcommand{\C}{{\mathbb C}}

\newcommand{\caT}{{\mathcal T}}
\newcommand{\tB}{{\tilde B}}

\DeclareMathOperator{\Res}{{Res}} \DeclareMathOperator{\wt}{{wt}}
\DeclareMathOperator{\sh}{{Shfl}}
\newcommand{\suf}{{\ast\,}}
\newcommand{\sufq}{{\ast_q\,}}
\newcommand{\gam}{{\gamma}}
\newcommand{\ga}{{\alpha}}

\newcommand{\gb}{{\beta}}
\newcommand{\gd}{{\delta}}
\newcommand{\gs}{{\sigma}}
\newcommand{\sif}{{\mathcal S}}
\newcommand{\gt}{{\tau}}

\newcommand{\fS}{{\mathfrak S}}
\newcommand{\DD}{{\mathfrak D}}

\newcommand{\zq}{{\zeta_q}}
\newcommand\qup{{q\uparrow 1}}

\begin{document}

\title{$q$-Multiple Zeta Functions and $q$-Multiple Polylogarithms}
\author{Jianqiang Zhao\footnote{Partially supported by
NSF grant DMS0139813}}
\date{}
\maketitle
\begin{center}
{\large Department of Mathematics, University of Pennsylvania, PA
19104, USA}
\end{center}

\noindent{\bf Abstract.}
We shall define the $q$-analogs of multiple zeta functions and
multiple polylogarithms in this paper and study their properties,
based on the work of Kaneko et al. and Schlesinger, respectively.

\section{Introduction and definitions}
Let $0<q<1$ and for any positive integer $k$ define its $q$-analog
$[k]=[k]_q=(1-q^k)/(1-q)$. In \cite{KNW} Kaneko et al. define a
function of two complex variables $f_q(s;t)=\sum_{k=1}^\infty
q^{kt}/[k]^s$ such that the $q$-analog of Riemann zeta function is
realized as
$$\zq(s):=f_q(s;s-1).$$
For any property $\mathcal P$ let $\Z_{\mathcal P}$ be the set of
integers satisfying $\mathcal P$. After analytically continuing
$f_q(s;t)$ to $\C^2$ as a meromorphic function Kaneko et al.
proved the following main result
\begin{thm}\label{thm:knw} {\em (\cite[Prop.~2, Thm.~2]{KNW})}
One can analytically continue $\zq(s)$ to $\C\setminus\{\Z_{\le
0}+\frac{2\pi i}{\log q}\Z_{\ne 0}\}$. Moreover, for any $s\ne 1$
one has
$$\lim_\qup \zq(s)=\zeta(s).$$
\end{thm}
They also study the special values of $\zeta_q(s)$ at non-negative
integers. In this paper we shall generalize these to the
(Euler-Zagier) multiple zeta functions, which are defined as
nested generalizations of Riemann zeta function $\zeta(s)$:
\begin{equation}\label{zeta}
\zeta(s_1,\dots, s_d)=\sum_{0<k_1<\dots<k_d} k_1^{-s_1}\cdots
k_d^{-s_d}
\end{equation}
for complex variables $s_1,\dots, s_d$ satisfying
$\gs_j+\dots+\gs_d>d-j+1$ for all $j=1,\dots,d$. Here and in what
follows, whenever $s\in \C$ we always write $\gs=\Re(s)$, the real
part of $s$. The analytic continuation of multiple zeta functions
has been studied independently in \cite{AET} and \cite{Zana}. We
know that $\zeta(s_1,\dots, s_d)$ can be extended to a meromorphic
function on $\C^d\setminus \fS_d$ where
\begin{equation}\label{fSd}
\fS_d=\left\{(s_1,\dots,s_d)\in \C^d\left| \aligned
s_d= &1,s_{d-1}+s_{d}=2,1,-2m \ \forall m\in\Z_{\ge 0},\\
\text{and }&s_j+\dots+s_d\in\Z_{\le d-j+1} \ \forall j\le d-2.
\endaligned
\right.\right\}
\end{equation}

To find the $q$-analog of multiple zeta functions we first define
an auxiliary function of $2d$ complex variables $s_1,\dots, s_d,
t_1,\dots, t_d\in \C$
\begin{equation*}
f_q(s_1,\dots, s_d;t_1,\dots, t_d)=\sum_{0<k_1<\dots<k_d}
\frac{q^{k_1t_1+\cdots+k_dt_d}}{[k_1]^{s_1}\cdots[k_d]^{s_d}}
\end{equation*}
which converges if $\Re(t_j+\dots+t_d)>0$ for all $j=1,\dots,d$
(see Prop.~\ref{prop:form}). In the next section we are going to
analytically continue this function to $\C^{2d}$ as a meromorphic
function with explicitly defined poles.

We now define the $q$-multiple zeta function by specialization of
$f_q$:
$$\zq(s_1,\dots,s_d):=f_q(s_1,\dots,s_d;s_1-1,\dots,s_d-1)$$
which will be shown to be the correct $q$-analog of multiple zeta
functions. When $\gs_j>1$ for all $j$ we can express this by the
series
\begin{equation}\label{qzeta}
\zq(s_1,\dots,s_d)=\sum_{0<k_1<\dots<k_d}
\frac{q^{k_1(s_1-1)+\cdots+k_d(s_d-1)}}{[k_1]^{s_1}\cdots[k_d]^{s_d}}.
\end{equation}
Note that when $d=1$ this is the same as the $q$-analog of the
Riemann zeta function defined in \cite{KNW}. Put
\begin{equation*}
\fS_d'=\left\{(s_1,\dots,s_d)\in \C^d: \aligned &s_d\in
1+\frac{2\pi i}{\log q}\Z, \text{ or }
s_d\in\Z_{\le 0}+\frac{2\pi i}{\log q}\Z_{\ne 0}, \\
&\text{or } s_j+\dots+s_d\in \Z_{\le d-j+1} +\frac{2\pi i}{\log
q}\Z,\  j<d
\endaligned\right\}\supset \fS_d.
\end{equation*}
Here the last part in $\fS_d'$ is vacuous if $d=1$. The primary
goal of this paper is to prove

\medskip
\noindent {\bf Main Theorem.} {\em The $q$-multiple zeta function
$\zq(s_1,\dots,s_d)$ can be extended to a meromorphic function on
$\C^d\setminus \fS'_d$ with simple poles along $\fS'_d$. Further,
for all} $(s_1,\dots,s_d)\in \C^d\setminus \fS_d$
\begin{equation*}
\lim_\qup \zq(s_1,\dots,s_d)=\zeta(s_1,\dots,s_d).
\end{equation*}

\medskip
In section \ref{sec:polylog}, we propose a new definition of the
$q$-multiple polylogarithms and briefly study their properties. We
also review Jackson's $q$-derivatives and $q$-definite integrals
and define $q$-iterated integrals as $q$-analogs of Chen's
iterated integrals.

It is known that there're two kinds of shuffle relations among
multiple zeta values (MZV for short). The first one is produced by
their power series expansions, the second by using Chen's iterated
integrals. In the last section of this paper we will apply our
$q$-iterated integral technique to $q$-multiple polylogarithms in
order to study the $q$-shuffle relations of the second kind for
$q$-MZV. For simplicity we will only deal with $\zq(m)\zq(n)$ for
positive integers $m\ne n$. These relations reduce to the ordinary
ones when $\qup$. We thank Prof. Kaneko for his questions relating
to this part of our study and sending us his offprint upon which
the current work is based.

\section{Analytic continuations of $f_q$ and $\zq$}\label{sec:fqzq}
The purpose of this section is two-fold: first we will use the
auxiliary functions $f_q$ introduced in the first section to give
a quick analytic continuation of $q$-multiple zeta functions
$\zq(s_1,\dots,s_d)$, though this is not enough to show it's the
right $q$-analog of the multiple zeta functions. Second, we write
down these expressions involving binomial coefficients explicitly
which will be used to study special values of $\zq(s_1,\dots,s_d)$
in section~\ref{sec:val}.

We need a simple lemma first.
\begin{lem}\label{lem:fracx}
Let $k$ be a positive integer. For all $1>q>1/2$ we have
$$1\le \frac{1-q^k}{1-q}< 2k.$$
\end{lem}
\begin{proof}
The inequality $1\le \frac{1-q^k}{1-q}$ is obvious. Let
$f(q)=2k(1-q)-(1-q^k)$. We show that $f(q)\ge 0$ for all $0<q<1$.
Indeed, for all such $q$
\begin{equation*}
f'(q)=-2k+kq^{k-1}\le -k<0.
\end{equation*}
So $f$ is strictly decreasing and the positivity of $f$ follows
from $f(1)=0$. This implies that $2k>(1-q^k)/(1-q)$ as desired.
\end{proof}

\begin{prop}\label{prop:form}
The function
\begin{equation*}
f_q(s_1,\dots, s_d;t_1,\dots, t_d)=\sum_{0<k_1<\dots<k_d}
\frac{q^{k_1t_1+\cdots+k_dt_d}}{[k_1]^{s_1}\cdots[k_d]^{s_d}}
\end{equation*}
converges if $\Re(t_j+\dots+t_d)>0$ for all $j=1,\dots,d$. It can
be analytically continued to a meromorphic function over $\C^{2d}$
via the series expansion
\begin{equation}\label{merom}
f_q(s_1,\dots, s_d;t_1,\dots, t_d)=(1-q)^{\wt(\ors)}
\sum_{r_1,\dots,r_d=0}^{+\infty} \prod_{j=1}^d
\left[{s_j+r_j-1\choose r_j}
\frac{q^{j(r_j+t_j)}}{1-q^{r_j+t_j+\cdots+r_d+t_d}} \right],
\end{equation}
where $\wt(\ors)=s_1+\dots+s_d$. It has the following (simple)
poles: $t_j+\dots+t_d\in \Z_{\le 0}+\frac{2\pi i}{\log q}\Z.$
\end{prop}
\begin{proof}
Assume $|\Re(s_j)|<N_j$ and let $\tau_j=\Re(t_j)$ for all
$j=1,\dots, d$. By Lemma~\ref{prop:form}
\begin{equation}\label{fbound}
\Bigl| f_q(s_1,\dots, s_d;t_1,\dots, t_d)\Bigr| <
\sum_{0<k_1<\dots<k_d} \prod_{j=1}^d (2k_j)^{N_1} q^{k_j\gt_j}.
\end{equation}
Let $k=k_{d-1}$ (when $d=1$ take $k_0=0$), $n=k_d$, $N=N_d$, and
$\gt=\gt_d$. Then by root test $\sum_{n>k}n^Nq^{n\gt}$ converges
and moreover
\begin{equation}\label{geomseries}
\sum_{n>k}n^Nq^{n\gt}=\left(q\frac{d}{dq}\right)^N
\sum_{n>k}q^{n\gt} =(1-q^\gt)^{-N}f_N(k;q,\gt) q^{k\gt}
\end{equation}
where $f_N(x;q,\gt)=\sum_{l=0}^N c_lx^l$ is a polynomial of degree
$N$ whose coefficients depend only on the constants $N$, $q$ and
$\gt$. This proves the first part of the lemma when $d=1$. In the
general case it follows from \eqref{fbound} and \eqref{geomseries}
that
\begin{equation*}
\Bigl| f_q(s_1,\dots, s_d;t_1,\dots, t_d)\Bigr| <
\frac{2^{N_1+\dots+N_d}}{(1-q^\gt)^N} \sum_{l=0}^{N} c_l
\sum_{0<k_1<\dots<k_{d-1}<k} \left(\prod_{j=1}^{d-2} (k_j)^{N_j}
q^{k_j\gt_j}\right) k^{N_{d-1}+l} q^{k(\gt_{d-1}+\gt_d)}.
\end{equation*}
Hence the first part of the lemma follows from an easy induction
on $d$.

By binomial expansion
$(1-x)^{-s}=\sum_{r=0}^{\infty}{s+r-1\choose r} x^r$ we get
 \begin{equation*}
f_q(s_1,\dots, s_d;t_1,\dots, t_d)=(1-q)^{\wt(\ors)}
\sum_{0<k_1<\dots<k_d} \sum_{r_1,\dots,r_d=0}^{+\infty}
\left(\prod_{j=1}^d {s_j+r_j-1\choose r_j}q^{k_j(r_j+t_j)}\right).
\end{equation*}
As $0<q<1$ the series converges absolutely by Stirling's formula
so we can exchange the summations. The proposition follows
immediately from the next lemma by taking $x_j=q^{t_j+r_j}$ for
$j=1,\dots,d$.
\end{proof}
\begin{lem}\label{lem:uuu} Let $x_j\in \C$ such
that $|x_j|<1$ for $j=1,\dots,d$. Then
\begin{equation}\label{equ:lem}
\sum_{0< k_1<\dots<k_d} \prod_{j=1}^d x_j^{k_j}=\prod_{j=1}^d
\frac{x_j^j }{1-x_j }=\prod_{j=1}^d \frac{x_j\cdots x_d
}{1-x_j\cdots x_d }.
\end{equation}
\end{lem}
\begin{proof}
When $n=1$ this is clear. For $n\ge 2$ we denote the right hand
side of \eqref{equ:lem} by $F_d(x_1,\dots,x_d)$. Then
\begin{equation*}
\sum_{0<k_1<\dots<k_d} \prod_{j=1}^d x_j^{k_j}=\sum_{0<
k_1<\dots<k_{d-1}} \prod_{j=1}^{d-1} x_j^{k_j} \sum_{k_d>k_{d-1}}
x_d^{k_d} =
 \frac{x_d}{1-x_d} F_{d-1}(x_1,\dots,x_{d-2},x_{d-1}x_d).
\end{equation*}
The lemma follows immediately by induction.
\end{proof}

Recall that
$\zq(s_1,\dots,s_d)=f_q(s_1,\dots,s_d;s_1-1,\dots,s_d-1)$. Hence
we have the following immediate consequence.
\begin{thm}\label{thm:qmz}
The $q$-multiple zeta function $\zq(s_1,\dots,s_d)$ can be
extended to a meromorphic function with simple poles lying along
$\fS_d'$:
\begin{equation*}
\zq(s_1,\dots, s_d)=(1-q)^{\wt(\ors)}
\sum_{r_1,\dots,r_d=0}^{+\infty} \prod_{j=1}^d
\left[{s_j+r_j-1\choose r_j}
\frac{q^{j(r_j+s_j-1)}}{1-q^{r_j+s_j+\cdots+r_d+s_d-d+j-1}}
\right].
\end{equation*}
\end{thm}

To see the effect of taking different specializations of $t_j$ in
$f_q$ we define the shifting operators $\sif_j$ ($1\le j\le d$) on
the multiple zeta functions as follows:
\begin{equation}\label{equ:shift}
\sif_j\zeta(s_1,\dots,s_d)=\zeta(s_1,\dots,s_d)+(1-q)\zeta(s_1,\dots,s_j-1,\dots,s_d).
\end{equation}
It is obvious that these operators are commutative.

\begin{prop} \label{shift}
Let $n_1,\dots,n_d$ be non-negative integers. Then we have
\begin{multline*}
f_q(s_1,\dots,s_d;s_1-1-n_1,\dots,s_d-1-n_d)\\
= \sif_1^{n_1}\circ\dots\circ\sif_d^{n_d} \zeta(s_1,\dots,s_d)=
\sum_{r_1=0}^{n_1}\cdots \sum_{r_d=0}^{n_d} \left(\prod_{j=1}^d
{n_j\choose r_j}(1-q)^{r_j}\right)\zq(s_1-r_1,\dots,s_d-r_d).
\end{multline*}
\end{prop}
\begin{proof} We only sketch the proof in the case $n_1=\cdots=n_{d-1}=0$.
The general case is completely similar. In the rest of the paper
we always let $\sif$ be the shifting operator on the last
variable. Suppose $n_d=n=1$. Then
\begin{align*}
\ &f_q(s_1,\dots,s_d;s_1-1,\dots,s_{d-1}-1,s_d-2)\\
=&\sum_{0<k_1<\dots<k_d}
\frac{q^{k_1(s_1-1)+\cdots+k_{d-1}(s_{d-1}-1)+k_d(s_d-2)}}{[k_1]^{s_1}\cdots[k_d]^{s_d}}\\
=& \sum_{0<k_1<\dots<k_d}
\frac{q^{k_1(s_1-1)+\cdots+k_{d-1}(s_{d-1}-1)}}
{[k_1]^{s_1}\cdots[k_{d-1}]^{s_{d-1}}}\cdot
\frac{q^{k_d(s_d-2)}(1-q^{k_d})+q^{k_d(s_d-1)}}{[k_d]^{s_d}}\\
=&\sif\zeta(s_1,\dots,s_d).
\end{align*}
The rest follows easily by induction.
\end{proof}
The next corollary answers an implicit question in
\cite{KNW}.
\begin{cor}\label{cor:shift}
Let $n$ be a positive integer. The specialization of $t$ in
$f_q(s;t)$ to $s-1-n$ is
\begin{equation*}
f_q(s;s-1-n)=\sif^n\zq(s)=\sum_{r=0}^n {n\choose r} (1-q)^r
\zq(s-r).
\end{equation*}
\end{cor}
We observe that one effect of the shifting operator is to bring in
more poles. Essentially, $\sif^n$ shifts all the poles of $\zq(s)$
by $n$ to the right on the complex plane.

\section{Analytic continuation of multiple zeta functions}
Let's begin with a review of some classical results on Bernoulli
polynomials $B_k(x)$ which are defined by the generating function
\begin{equation*}
\frac{te^{xt}}{e^t-1}=\sum_{k=0}^\infty B_k(x)\frac{t^k}{k!}.
\end{equation*}
Let $\tB_k(x)$ be the ``periodic Bernoulli polynomial''
\begin{equation*}
\tB_k(x)=B_k(\{x\}), \quad x\ge 1,
\end{equation*}
where $\{x\}$ is the fractional part of $x$. Then we have
(\cite[Ch.~IX, Misc.~Ex.~12]{WW})
\begin{equation} \label{equ:perber}
\tB_k(x)=-k!\sum_{n\in \Z\setminus\{0\}} \frac{e^{2\pi inx}}{(2\pi
in)^k}.
\end{equation}
Recall that the Bernoulli numbers satisfy $B_k=B_k(1)$ if $k\ge 2$
while $B_0=1$ and $B_1=-1/2$.

\begin{lem} \label{lem:berbound}
For every positive integer $M\ge 2$ and $x>1$ we have
$$\Bigl|\tB_M(x)\Bigr| \le \frac{4 M!}{(2\pi)^M}.$$
\end{lem}
\begin{proof}
It follows from the fact that $\zeta(M)\le \zeta(2)=\pi^2/6<2$ for
$M\ge 2$.
\end{proof}

We know that one can analytically continue the multiple zeta
functions as independently presented in \cite{Zana} and \cite{AET}
by different methods. Moreover, $\zeta(s_1,\dots, s_d)$ has
singularities on the hyperplanes in $\fS_d$ defined by
\eqref{fSd}. However, neither approach is suitable for our purpose
here. So we follow the idea in \cite{KNW} to provide a third
approach in the rest of this section. The same idea will also be
used to deal with the $q$-multiple zeta functions.

Let's recall the classical Euler-Maclaurin summation formula
\cite[7.21]{WW}. Let $f(x)$ be any (complex-valued) $C^\infty$
function on $[1,\infty)$ and let $m$ and $M$ be two positive
integers. Then we have
\begin{align}
\sum_{n=1}^m f(n)=&\int_1^m f(x)\,
dx+\frac{1}{2}(f(1)+f(m))+\sum_{r=1}^M
\frac{B_{r+1}}{(r+1)!} \left(f^{(r)}(m)-f^{(r)}(1)\right)\nonumber\\
-&\frac{(-1)^{M+1}}{(M+1)!}\int_1^m \tB_{M+1}(x)
f^{(M+1)}(x)\,dx.\label{equ:EM}
\end{align}

To simplify our notation, in definition \eqref{zeta} we replace
$s_d$, $k_{d-1}$, $k_d$ by $s$, $k$ and $n$, respectively. Taking
$f(x)=1/x^s$ and $m=k$ and $\infty$ in \eqref{equ:EM} we have:
\begin{align*}
\sum_{n>k} \frac{1}{n^s}=& \sum_{n=1}^\infty
\frac{1}{n^s}-\sum_{n=1}^{k} \frac{1}{n^s}\\
=&\int_{k}^\infty f(x)\, dx-\frac{1}{2}f(k)- \sum_{r=1}^M
\frac{B_{r+1}}{(r+1)!}f^{(r)}(k)-\frac{(-1)^{M+1}}{(M+1)!}\int_k^\infty
\tB_{M+1}(x)f^{(M+1)}(x)\,dx \\
=&\frac{1}{(s-1)k^{s-1}}-\frac{1}{2k^s}+ \sum_{r=1}^M
\frac{B_{r+1}}{(r+1)!}\frac{(s)_r}{k^{s+r}} -
\frac{(s)_{M+1}}{(M+1)!} \int_k^\infty
 \frac{\tB_{M+1}(x)}{x^{s+M+1}}\,dx.
\end{align*}
Here we have used the fact that $B_k=0$ if $k\ge 3$ is odd. By
definition \eqref{zeta} we have
\begin{thm}\label{thm:mz}
For all $(s_1,\dots,s_d)\in\C^d\setminus \fS_d$ and $M
>1+|\gs_d|+|\gs_{d-1}|$ we have
\begin{multline} \label{equ:key}
\zeta(s_1,\dots,s_d)=\sum_{r=0}^{M+1}
\frac{B_r}{r!} (s_d)_{r-1}\cdot \zeta(s_1,\dots,s_{d-1}+s_d+r-1)\\
- \frac{(s_d)_{M+1}}{(M+1)!} \sum_{0<k_1<\dots<k_d}
\frac{1}{k_1^{s_1}\cdots k_{d-1}^{s_{d-1}}} \int_{k_{d-1}}^\infty
\frac{\tB_{M+1}(x)}{x^{s_d+M+1}}\,dx.
\end{multline}
where we set $(s)_0=1$ and $(s)_{-1}=1/(s-1)$. This provides an
analytic continuation of $\zeta(s_1,\dots, s_d)$ to $\C^d\setminus
\fS_d$.
\end{thm}
\begin{proof} We only need to show that the series in
\eqref{equ:key} converges. Lemma \ref{lem:berbound} implies
(if $d=2$ then take $k_0=1$)
\begin{equation*}
\sum_{k_{d-1}=k_{d-2}}^\infty \left|\frac{1}{k_{d-1}^{s_{d-1}}}
\int_{k_{d-1}}^\infty \frac{\tB_M(x)}{x^{M+s_d+1}}\,dx\right|\le
\frac{4M!}{(2\pi)^M(M-|\gs_d|)} \sum_{k_{d-1}=k_{d-2}}^\infty
\frac{1}{k_{d-1}^{M-|\gs_{d-1}|-|\gs_d|}}
\end{equation*}
which converges absolutely whenever $M>1+|\gs_d|+|\gs_{d-1}|$.
\end{proof}

\section{Proof of Main Theorem}
Fix $(s_1,\dots,s_d)\in \C^d$ such that $\gs_j+\cdots+\gs_d<d-j+1$
for all $j=1,\dots,d.$ When $d=1$ this is Thm.~\ref{thm:knw}  due
to Kaneko et al. We now assume $d\ge 2$ and proceed by induction.
The key is a recursive formula for $\zq(s_1,\dots,s_d)$ similar to
\eqref{equ:key} for $\zeta(s_1,\dots,s_d)$. To derive this formula
we appeal to the Euler-Maclaurin summation formula \eqref{equ:EM}
again. Hence we set
\begin{equation*}
F(x)= \frac{q^{x(s-1)}}{(1-q^x)^s}
\end{equation*}
as in \cite{KNW}. Then
\begin{align*}
F'(x)=&(\log q)q^{x(s-1)} \frac{s-1+q^x}{(1-q^x)^{s+1}},\\
F''(x)=&(\log q)^2q^{x(s-1)}
\frac{s(s+1)-3s(1-q^x)+(1-q^x)^2}{(1-q^x)^{s+2}}.
\end{align*}

In definition \eqref{qzeta} we replace $s_d$, $k_{d-1}$, and $k_d$
by  $s$, $k$ and $n$, respectively. We now take $M=1$,
$f(x)=F(x+k-1)$  and let $m\to \infty$ in \eqref{equ:EM} and get
\begin{align}
\sum_{n>k} \frac{q^{n(s-1)}}{[n]^s}=&(1-q)^s\left(-F(k)+\sum_{n=1}^\infty f(n)\right) \nonumber \\
=&(1-q)^s\left(\int_{k}^\infty F(x)\, dx-\frac{1}{2}F(k)-
\frac{1}{12}F'(k)-\frac{1}{2}\int_1^\infty\tB_2(x)f''(x)\,dx\right) \nonumber \\
=&\frac{(q-1)}{(s-1)\log q}\frac{q^{k(s-1)}}{[k]^{s-1}}
 -\frac{1}{2}\frac{q^{k(s-1)}}{[k]^s}
 +\frac{1}{12}\frac{\log q}{q-1}
 \frac{q^{k(s-1)}(s+q^k-1)}{[k]^{s+1}} \nonumber \\
 -&\frac{(1-q)^s(\log q)^2}{2}\int_k^\infty\tB_2(x) q^{x(s-1)}
\frac{s(s+1)-3s(1-q^x)+(1-q^x)^2}{(1-q^x)^{s+2}}\,dx
\label{innersum}
\end{align}
because $\tB_2(x+k-1)=\tB_2(x)$ by periodicity. By the same
argument as in \cite{KNW}, setting the incomplete beta integrals
\begin{equation*}
b_t(\ga,\gb)=\int_0^t u^{\ga-1}(1-u)^{\gb-1} \, du,
\end{equation*}
we can obtain from \eqref{equ:perber} the following expression for
the last term in \eqref{innersum}:
\begin{equation}\label{lastterm}
\eqref{innersum}=
    -\sum_{n\in\Z\setminus\{0\}} \frac{(1-q)^s\log q}{(2\pi in)^2}
    \sum_{\nu=\pm 1,0}a_\nu(s) b_{q^k}(s-1+\gd n,-s+\nu)
\end{equation}
where $\gd=2\pi i/\log q$, $a_{-1}(s)=s(s+1),a_0(s)=-3s,$ and
$a_1(s)=1$. Repeatedly applying integration by parts on these
incomplete beta integrals we get for $\nu=\pm 1,0$ and positive
integer $M\ge 2$
\begin{multline}\label{bqk}
b_{q^k}(s-1+\gd n,-s+\nu)=\sum_{r=1}^{M-1} (-1)^{r-1}
\frac{(s+1-\nu)_{r-1}}{(s-1+\gd n)_r}
\frac{q^{k(s+r-2)}}{(1-q^k)^{s+r-\nu}} \\
+(-1)^{M-1} \frac{(s+1-\nu)_{M-1}}{(s-1+\gd n)_{M-1}}
b_{q^k}(s-2+M+\gd n, -s-M+1+\nu).
\end{multline}

Set $\ors'=(s_1,\dots,s_{d-2})$ if $d\ge 3$ and $\ors'=\emptyset$
if $d=2$. Putting \eqref{qzeta}, \eqref{innersum},
\eqref{lastterm} and \eqref{bqk} together and applying
Prop.~\ref{shift} we get
\begin{multline}\label{qzkey}
\zq(\ors',s_{d-1},s_d)=\frac{(q-1)}{(s-1)\log
q}\zq(\ors',s_{d-1}+s_d-1)
-\frac{1}{2} \sif\zq(\ors',s_{d-1}+s_d) \\
+\frac{s}{12}\frac{\log q}{q-1}\sif^2\zq(\ors',s_{d-1}+s_d+1)
+\frac{\log q}{12} \sif\zq(\ors',s_{d-1}+s_d) -\sum_{\nu=\pm
1,0}(C_\nu+D_\nu)
\end{multline}
where $C_\nu$ and $D_\nu$ are the contributions from the sum
involving $b_{q^k}(\dots,-s+\nu)$. Explicitly they are computed as
follows. Write
\begin{equation}\label{T}
T(q,s,n,r)=\prod_{j=0}^{r-1} \big(2\pi in+(s-1+j)\log q\big)^{-1}.
\end{equation}
Then
\begin{align*}
C_{-1}=&\sum_{r=1}^{M-1}  \sum_{n\in\Z\setminus\{0\}}
\frac{T(q,s,n,r)}{(2\pi i n)^2} \left(\frac{\log
q}{q-1}\right)^{r+1} (s)_{r+1}
\cdot \sif^3\zq(\ors',s_{d-1}+s_d+r+1), \\
D_{-1}=&-\sum_{0<k_1<\dots<k_{d-2}<k_{d-1}}
\frac{q^{k_1(s_1-1)+\cdots+k_{d-1}(s_{d-1}-1)}}
{[k_1]^{s_1}\cdots[k_{d-1}]^{s_{d-1}}} \sum_{n\in\Z\setminus\{0\}}
 \frac{T(q,s,n,M-1)}{(2\pi i n)^2} \cdot\\
\ &\qquad \cdot \left(\frac{\log q}{q-1}\right)^{M+1}(s)_{M+1}
 \cdot\int_{k_{d-1}}^\infty e^{2\pi in x} q^{x(s-2+M)}
    \left(\frac{1-q^x}{1-q}\right)^{-s-M-1}\,dx\\
   =&\sum_{0<k_1<\dots<k_{d-2}}
  \frac{q^{k_1(s_1-1)+\cdots+k_{d-2}(s_{d-2}-1)}}
{[k_1]^{s_1}\cdots[k_{d-2}]^{s_{d-2}}} \sum_{k=k_{d-2}+1}^\infty
R(M,q,k,s_{d-1},s).
\end{align*}
where we replace the index $k_{d-1}$ by $k$. Similarly,
\begin{align*}
C_0=&3\log q\sum_{r=1}^{M-1} \sum_{n\in\Z\setminus\{0\}}
 \frac{T(q,s,n,r)}{(2\pi i n)^2}\left(\frac{\log q}{q-1}\right)^r
 (s)_r\cdot\sif^2\zq(\ors',s_{d-1}+s_d+r),\\
D_0=&-3\log q\sum_{0<k_1<\dots<k_{d-2}<k_{d-1}}
\frac{q^{k_1(s_1-1)+\cdots+k_{d-1}(s_{d-1}-1)}}
{[k_1]^{s_1}\cdots[k_{d-1}]^{s_{d-1}}} \sum_{n\in\Z\setminus\{0\}}
 \frac{T(q,s,n,M-1)}{(2\pi i n)^2} \cdot\\
\ & \cdot\left(\frac{\log q}{q-1}\right)^M (s)_M\cdot
\int_{k_{d-1}}^\infty e^{2\pi in x}
 q^{x(s-2+M)} \left(\frac{1-q^x}{1-q}\right)^{-s-M}\,dx,
\end{align*}
and
\begin{align*}
C_1=&(\log q)^2 \sum_{r=1}^{M-1} \sum_{n\in\Z\setminus\{0\}}
 \frac{T(q,s,n,r)}{(2\pi i n)^2}\left(\frac{\log q}{q-1}\right)^{r-1}
 (s)_{r-1} \cdot \sif\zq(\ors',s_{d-1}+s_d+r-1), \\
D_1=&-(\log q)^2\sum_{0<k_1<\dots<k_{d-2}<k_{d-1}}
\frac{q^{k_1(s_1-1)+\cdots+k_{d-1}(s_{d-1}-1)}}
{[k_1]^{s_1}\cdots[k_{d-1}]^{s_{d-1}}} \sum_{n\in\Z\setminus\{0\}}
\frac{T(q,s,n,M-1)}{(2\pi i n)^2}\cdot \\
\ & \cdot  \left(\frac{\log q}{q-1}\right)^{M-1} (s)_{M-1}\cdot
\int_{k_{d-1}}^\infty e^{2\pi in x}
 q^{x(s-2+M)} \left(\frac{1-q^x}{1-q}\right)^{-s-M+1} \,dx.
\end{align*}

The crucial step next is to control the summations over $k_{d-1}$
and show that they converge uniformly with respect to $q$. When
$0<q\le 1/2$ this is clear. The only non-trivial part is when
$\qup$. So we assume $1/2<q<1$. Note that
$$\lim_\qup T(q,s,n,r)=\frac{1}{(2\pi i n)^r},\quad
\lim_\qup \frac{\log q}{q-1}=1.$$

\begin{lem}\label{inequ}
Let $s_d=s=\gs+i\tau$. Let $q_0=\max\{1/2, e^{(6-2\pi)/\tau}\}$ if
$\tau>0$ and let $q_0=1/2$ if $\tau\le 0$. Then for all $1>q>q_0$
and positive integer $k$ we have
$$\left|\frac{\log q}{q-1}\right|<2,
\quad\text{ and }\quad |T(q,s,n,r)|<\frac{1}{(6n)^r}.$$
\end{lem}
\begin{proof}
Let $f(q)=2(1-q)+\log q$. Then $f'(q)=-2+1/q<0$ whenever $q>q_0$.
So $f(q)>f(1)=0$ whenever $1>q>q_0$. This implies that
$2(1-q)>-\log q$ whence $\log q/(q-1)< 2$.

To bound $T(q,s,n,r)$ we consider each of its factors in
definition \eqref{T}. For each $0\le j<r$ we have
$$|2\pi in+(s-1+j)\log q|^2=\big((\gs-1+j)\log q\big)^2
 +(2\pi n+\tau\log q)^2\ge (2\pi n+\tau\log q)^2$$
which is independent $j$. If $\tau\le 0$ then clearly $|2\pi
in+(s-1+j)\log q|>6n$. If $\tau>0$ then it follows from
$q>e^{(6-2\pi)\tau}$ that
\begin{equation*}
    2\pi n+\tau\log q>2\pi n+6-2\pi\ge 6n,
\end{equation*}
as desired.
\end{proof}
Next we want to bound the integral terms in $D_{-1}$. Let
$|\gs_d|<N$ and $|\gs_{d-1}|<N'$ for some positive integers $N$
and $N'$. Fix an arbitrary $x>k$ and a positive integer
$M>16+2N+6\sum_{j=1}^{d-1}(N_j+1)$.  Then
\begin{equation*}
q^{-k(M/6-N'-1)} \left|q^{k(s'-1)}
 q^{x(s-2+M)} \left(\frac{1-q^x}{1-q}\right)^{-s-M-1}\right|
 < q^{x(M-N-2)-kM/6}
 \left(\frac{1-q}{1-q^x}\right)^{M-N+1}.
\end{equation*}
Denote by $g(q)$ the right hand side of the above inequality.
\begin{lem} \label{lem:lim}
Let $1/2<q<1$. Then $g(q)$ is increasing as a function of $q$ so
that
\begin{equation*}
    g(q)\le \lim_\qup g(q)=\frac{1}{x^{M-N+1}}.
\end{equation*}
\end{lem}
\begin{proof} Taking the logarithmic derivative of $g(q)$ we have
\begin{align*}
\frac{g'(q)}{g(q)}=&\frac{x(M-N-2)-kM/6}{q}+
(M-N+1)\frac{xq^{x-1}(1-q)-(1-q^x)}{(1-q)(1-q^x)}\\
=&\frac{(1-q)(1-q^x)(x(M-N-2)-kM/6)+(M-N+1)(xq^x(1-q)-q+q^{x+1})}
{q(1-q)(1-q^x)}
\end{align*}
whose numerator is denoted by $h(q)$. Then
\begin{align*}
 h'(q)=&\big((x+1)q^x-xq^{x-1}-1\big)\big(x(M-N-2)-kM/6)\big)\\
\ &+(M-N+1)\Bigl(x\big(xq^{x-1}-(x+1)q^x\big)-1+(x+1)q^x\Bigr).
\end{align*}
Clearly $h'(1)=0$ and moreover
\begin{alignat*}{3}
 q^{2-x}h''(q)=&\big(x(x+1)q-x(x-1)\big)\big(x(M-N-2)-kM/6)\big)\\
\ &+ (M-N+1)\Bigl(x\big(x(x-1)-x(x+1)q\big)+x(x+1)q\Bigr)\\
=&x(x-1)(kM/6+3x)+qx(x+1)\big(M-N+1-(kM/6+3x)\big)\\
=&(kM/6+3x)\big(x^2(1-q)-x(1+q)\big)+qx(x+1)(M-N+1) \\
>&qx(x+1)(M-N+1)-x(1+q)(kM/6+3x) &(\text{since }1>q)\\
\ge & qx^2\Bigl\{M-N+1-\frac{1+q}{q}(3+M/6)\Bigr\} &(\text{since }k<x)\\
> & qx^2(M/2-N-8)>0
\end{alignat*}
where we used the fact that if $q>1/2$ then $(1+q)/q<3$. This
implies that $h'(q)$ is increasing so that $h'(q)<0$ for all
$1>q>1/2$ (recall that $h'(1)=0$). It follows that $h(q)$ is
decreasing. But $h(1)=0$ so we know $h(q)>0$ for all such $q$.
Thus $g'(q)>0$ and therefore $g(q)$ is increasing. This completes
the proof of the lemma.
\end{proof}

We now can bound the innermost sum of $D_{-1}$. From
Lemma~\ref{lem:fracx}, Lemma~\ref{inequ} and Lemma~\ref{lem:lim}
we have (if $d=2$ then take $k_0=1$)
\begin{align*}
\sum_{k=1+k_{d-2}}^\infty | R(M,q,k,s',s) |
< & \sum_{k=1+k_{d-2}}^\infty  (2k)^{N'} \frac{2\zeta(M+1)}{4\pi^2
6^{M-1}} 2^{M+1} (N)_{M+1} q^{k(M/6-N'-1)}
\int_k^\infty \frac{dx}{x^{M-N+1}}\\
< & \frac{(M+1)!}{M-N} {M+N\choose M+1}
\sum_{k=1+k_{d-2}}^\infty \frac{q^{k(M/6-N'-1)}}{k^{M-N-N'}}  \\
<&(M+1)! (2M)^{M+1}  \sum_{k=1+k_{d-2}}^\infty q^{k(M/6-N'-1)}
\end{align*}
since $2^{N'+M+2}\zeta(M+1)<4\pi^2 6^{M-1}$ and $M-N>2$. Therefore
by Lemma~\ref{lem:fracx}
\begin{equation*}
\big|D_{-1}\big|<  (2M)^{2M} \sum_{0<k_1<\dots<k_{d-2}}
\left(\prod_{l=1}^{d-2}k_l^{N_l}q^{k_l(-N_l-1)}\right)
q^{k_{d-2}(M/6-N_{d-1}-1)}
\end{equation*}
which converges as proved in Prop.~\ref{prop:form}.

Exactly the same argument applies to the integral terms in $D_0$
and $D_1$, which we leave to the interested readers. These
convergence results imply two things. First we can show by
induction on $d$ that \eqref{qzkey} gives rise to an analytic
continuation of $\zq(s_1,\dots,s_d)$ as a meromorphic function on
$\C^d\setminus \fS_d'$. Second, also by induction on $d$, we now
can conclude that it's legitimate to take the limit $\qup$ inside
the sums of $C_\nu$ and $D_\nu$ to get (note that $\lim_\qup
\sif^n \zq(\ors)=\zeta(\ors)$ for any
$\ors\in\C^{d-1}\setminus\fS_{d-1}'$ and any positive integer $n$)
\begin{align*}
\lim_\qup\zq(\ors',s_{d-1},s_d)=&\frac{1}{s-1}\zeta(\ors',s_{d-1}+s_d-1)
-\frac{1}{2}\zeta(\ors',s_{d-1}+s_d)
+\frac{s}{12}\zeta(\ors',s_{d-1}+s_d+1) \\
\ &-\sum_{r=1}^{M-1} \sum_{n\in\Z\setminus\{0\}}
\frac{1}{(2\pi i n)^{r+2}} (s)_{r+1} \cdot\zeta(\ors',s_{d-1}+s_d+r+1) \\
\ &+\sum_{0<k_1<\dots<k_{d-1}} \frac{1}{k_1^{s_1}\cdots
k_{d-1}^{s_{d-1}}} \sum_{n\in\Z\setminus\{0\}}
 \frac{(s)_{M+1}}{(2\pi i n)^{M+1}}
\int_{k_{d-1}}^\infty e^{2\pi in x} x^{-s-M-1}\,dx\\
=&\sum_{r=0}^{M+1}
\frac{B_r}{r!} (s_d)_{r-1}\cdot \zeta(s_1,\dots,s_{d-1}+s_d+r-1)\\
\ &- \frac{1}{(M+1)!} \sum_{0<k_1<\dots<k_d}
\frac{1}{k_1^{s_1}\cdots k_{d-1}^{s_{d-1}}} \int_{k_{d-1}}^\infty
\tB_{M+1}(x)
 \frac{(s)_{M+1}}{x^{s+M+1}}\,dx
\end{align*}
by \eqref{equ:perber} and its specialization with $x=1$
$$\sum_{n\in \Z\setminus\{0\}}
\frac{1}{(2\pi in)^{r+2}}=-\frac{\tB_{r+2}(1)}{(r+2)!}
=-\frac{B_{r+2}}{(r+2)!}.$$ The main theorem now mostly follows
from Thm.~\ref{thm:mz}. The poles at $s_d=m-\frac{2\pi i}{\log
q}n$ are given by the first term in formula \eqref{qzkey} when
$m=1$ and $n=0$ and by the terms $T(q,s_d,n,r)$ as defined in
\eqref{T} if $m\le 1$ and $n\ne 0$. The location of the other
poles are obtained by induction using those poles of the
$q$-Riemann zeta function presented in Thm.~\ref{thm:knw} for the
initial step. This completes the proof of our main theorem.

\section{Series $q$-shuffle relations}\label{sec:shuffle1}
The classical multiple zeta functions satisfy shuffle relations
originating from their series representations. For example,
\begin{equation}\label{shufdb}
\zeta(s_1)\zeta(s_2)=\zeta(s_1,s_2)+\zeta(s_2,s_1)+\zeta(s_1+s_2).
\end{equation}
In general we can define a shuffle operation $\suf$ on finite
words so that for any complex numbers $a$, $b$ ( regarded as
letters) and words $w_1$ and $w_2$ of complex numbers
\begin{equation*}
aw_1\suf bw_2=a(w_1\suf bw_2)+b(aw_1\suf w_2)+(a+b)(w_1\suf w_2)
\end{equation*}
so that treating complex variables as words we have
$$\zeta(w_1)\zeta(w_2)=\zeta(w_1\suf w_2).$$

To generalize these relations we first define shifting operators
on words of complex variables:
\begin{equation*}
(s_1,\dots,s_{i-1},\sif(s_i),s_{i+1},\dots,s_d)=\sif_i
(s_1,\dots,s_d)
\end{equation*}
and set $\zeta\big(\sif_i (s_1,\dots,s_d)\big)=\sif_i
\zeta(s_1,\dots,s_d)$. Now we define the $q$-shuffle operator
$\sufq$ on words by
\begin{equation*}
aw_1\sufq bw_2=a(w_1\sufq bw_2)+b(aw_1\sufq
w_2)+\sif(a+b)(w_1\sufq w_2).
\end{equation*}
\begin{thm} For any two words of complex variables $w_1$ and $w_2$
we have
\begin{equation*}
\zq(w_1)\zq(w_2)=\zq(w_1 \sufq w_2).
\end{equation*}
\end{thm}
\begin{proof} Induction on the length of $w_1w_2.$
\end{proof}

For example
\begin{equation}\label{equ:dbshf}
\zq(s_1)\zq(s_2)
=\zq(s_1,s_2)+\zq(s_2,s_1)+\zq(s_1+s_2)+(q-1)\zq(s_1+s_2-1).
\end{equation}
We can recover the shuffle relation of \eqref{shufdb} by taking
$\lim_\qup$.

\section{Special values of $\zq(s_1,s_2)$}\label{sec:val}
For integers $n_1,\dots,n_d$ we set
\begin{align*}
\zq(n_1,\dots, n_d)=\lim_{s_1\to n_1}\cdots
\lim_{s_d\to n_d}\zq(s_1,\dots,s_d),\\
\zeta^R_q(n_1,\dots, n_d)=\lim_{s_d\to n_d}\cdots \lim_{s_1\to
n_1}\zq(s_1,\dots,s_d)
\end{align*}
if the limits exist. Interesting phenomena occur already in the
case $d=2$ and these should be generalized to arbitrary depth. By
Thm.~\ref{thm:qmz} we get
\begin{align*}
\zq(s_1, s_2)=&(1-q)^{s_1+s_2} \sum_{r_1,r_2=0}^{+\infty}
 {s_1+r_1-1\choose r_1} {s_2+r_2-1\choose r_2}
\frac{q^{2s_2+2r_2+s_1+r_1-3}}{(1-q^{s_2+r_2-1})(1-q^{s_2+r_2+s_1+r_1-2})}\\
=&(1-q)^{s_1+s_2} \left\{
\frac{q^{2s_2+s_1-3}}{(1-q^{s_2-1})(1-q^{s_2+s_1-2})}
+\frac{s_1q^{2s_2+s_1-2}}{(1-q^{s_2-1})(1-q^{s_2+s_1-1})}.\right.\\
+&\frac{s_2q^{2s_2+s_1-1}}{(1-q^{s_2})(1-q^{s_2+s_1-1})}
+\frac{s_1s_2q^{2s_2+s_1}}{(1-q^{s_2})(1-q^{s_2+s_1})} \\
+&\left.\frac{s_1(s_1+1)q^{2s_2+s_1-1}}{2(1-q^{s_2-1})(1-q^{s_2+s_1})}
+\frac{s_2(s_2+1)q^{2s_2+s_1+1}}{2(1-q^{s_2+1})(1-q^{s_2+s_1})}
+\cdots \right\}
\end{align*}
Clearly we have
\begin{align*}
\zq(0,0)=\lim_{s_1\to 0}\lim_{s_2\to 0}\zq(s_1,
s_2)=&\frac{1}{(q^2-1)(q-1)}-\frac{3}{2(q-1)\log
q}+\frac{1}{\log^2 q},\\
\zq^R(0,0)=\lim_{s_2\to 0}\lim_{s_1\to 0}\zq(s_1,
s_2)=&\frac{1}{(q^2-1)(q-1)}-\frac{1}{(q-1)\log q}+
\frac{q}{2(q-1)\log q}.
\end{align*}
It is not too hard to find that
$$
\lim_\qup\zq(0,0)=\frac{1}{3},\quad
\lim_\qup\zq^R(0,0)=\frac{5}{12}.
$$
This is consistent with what we found in \cite{Zana} by using
generalized functions (distributions). See also equations
\eqref{zeta00} below. In \cite{Zana} we further showed that near
$(0,0)$ the double zeta function has the following asymptotic
expansion:
$$\zeta(s_1,s_2)=\frac{4s_1+5s_2}{12(s_1+s_2)}+R(s_1,s_2)$$
where $R(s_1,s_2)$ is analytic at $(0,0)$ and
$\lim_{(s_1,s_2)\to(0,0)} R(s_1,s_2)=0.$

Let $n$, $k$ be two non-negative integers, and $m=k-n-2$. We now
consider the double zeta function around $(s_1,s_2)=(-m,-n)$ which
has the following expression by Thm.~\ref{thm:mz}:
\begin{multline} \label{equ:key1}
\frac{-1}{n+1}\zeta(s_1+s_2-1)-\frac{1}{2}\zeta(s_1+s_2)
-\sum_{r=2}^{n+1} \frac{B_r}{r} {n\choose r-1} \zeta(s_1+s_2+r-1)\\
+\ga(m)(-1)^n \frac{B_k}{k!} n!(k-n-2)!
(s_2+n)\zeta(s_2+s_1+k-1)
\end{multline}
where $\ga(m)=0$ if $m\le -1$ and $\ga(m)=1$ if $m\ge 0$.
Note that the last term is zero when computing
\begin{equation*}
\zeta(-m,-n)=\lim_{s_1\to -m}\ \lim_{s_2\to -n}\zeta(s_1,s_2)
\end{equation*}
while it has possibly nontrivial contribution for $\zeta^R(-m,-n)$ since
\begin{equation*}
\lim_{s_2\to -n} \lim_{s_1\to -m} (s_2+n)\zeta(s_2+s_1+k-1)=1.
\end{equation*}
We get
\begin{table}[h]
\begin{center}
\begin{tabular}{  ||c|c |c|| }
 \hline
 $k,\ n$& pole,residue & indeterminacy,
 $\zeta=\zeta^R(n+2-k,-n)$\\ \hline
 $k=0$ & $-1/(n+1)$ & NO \\  \hline
 $k=1$ & $-1/2$ & NO \\  \hline
 $2\nmid k, 3\le k\le n+1$ & NO & $\displaystyle{
  \frac{B_{k-1}}{2(k-1)}+\sum_{r=k-1}^{n+1}
     \frac{B_r}{r} {n\choose r-1}\zeta(r+1-k)}$ \\ \hline
 $2\nmid k, k= n+2$ & NO & $B_{k-1}/(k-1)$ \\  \hline
 $2\nmid k, k> n+2$ & NO & $B_{k-1}/2(k-1)$ \\  \hline
 $2|k, 2\le k\le n+1$ & $\displaystyle{(-1)^{k+1}\frac{B_{k}}{k} {n\choose k-1}} $ &NO \\ \hline
\end{tabular}
\caption{Poles and indeterminacy of double zeta function.}
\label{Ta:dbzeta}
\end{center}
\end{table}

When $m\ge 0$ and $2|k$ the values of $\zeta(-m,-n)$ and
$\zeta^R(-m,-n)$ are different in general:
\begin{align}
\zeta(-m,-n)=& \frac{B_{k}}{k(n+1)}+\sum_{r=1}^{n+1}
  \frac{B_r}{r} {n\choose r-1}\frac{B_{k-r}}{k-r}\label{equ:kn+2}\\
\zeta^R(-m,-n)=& \zeta(-m,-n)+(-1)^n \frac{B_k}{k!}
n!(k-n-2)!.\label{equ:kn+2R}
\end{align}
Note that the term corresponding to $r=1$ is non-zero if and only
if $k=2$ (and $n=0$). From this observation we again recover that
\begin{equation}\label{zeta00}
\zeta(0,0)=\frac{1}{3}, \quad \zeta^R(0,0)=\frac{5}{12}.
\end{equation}

We now consider the $q$-double zeta function.
\begin{thm} \label{thm:kn+1}
Let $k$, $n$ be two non-negative integers, and $m=k-n-2$.
If $m\le -1$ then the $q$-double zeta function $\zq(s_1,s_2)$
has a pole at $(-m,-n)$ with residue given by:
\begin{equation*}
\frac{\underset{(s_1,s_2)=(-m,-n)}{\Res}
\zq(s_1,s_2)}{-(1-q)^{2-k}(\log q)^{-1}}=
\begin{cases}
\displaystyle{\sum_{r=0}^k (-1)^r {n+1-r\choose k-r}{n\choose
r}/(q^{n+1-r}-1) }&
\quad\text{ if }k\le n, \\
\displaystyle{\sum_{r=0}^n (-1)^r{n\choose
r}/(q^{n+1-r}-1)-\frac{(-1)^n}{(n+1)\log q}} & \quad\text{ if }k=
n+1.
\end{cases}
\end{equation*}
\end{thm}
\begin{proof}
Use Thm.~\ref{thm:qmz}.
\end{proof}
\begin{cor} \label{cor:compare1}
Let $n$ be a non-negative integer. Then
\begin{equation}\label{equ:zn}
\underset{(s_1,s_2)=(1,-n)}{\Res} \zq(s_1,s_2)=\frac{q-1}{\log
q}\zq(-n)
\end{equation}
and
\begin{equation}\label{equ:zn1}
\lim_\qup \zq(-n)=-\frac{B_{n+1}}{n+1}=
\underset{(s_1,s_2)=(1,-n)}{\Res}\zeta(s_1,s_2).
\end{equation}
\end{cor}
\begin{proof} Equation \eqref{equ:zn} follows from
the case $m=-1$ in the above theorem and \cite[(6)]{KNW}:
\begin{equation*}
\zq(-n)=(1-q)^{-n}\left\{\sum_{r=0}^n (-1)^r{n\choose
r}/(q^{n+1-r}-1)-\frac{(-1)^n}{(n+1)\log q} \right\}.
\end{equation*}
The first equality in \eqref{equ:zn1} is \cite[Thm.~1]{KNW} and
the second equality follows from Table~\ref{Ta:dbzeta}.
\end{proof}

\begin{cor} \label{cor:compare2}
Let $k$, $n$ be two non-negative integers, $m=k-n-2\le -2$. Then
\begin{equation}\label{equ:zkleqn}
\lim_\qup \underset{(s_1,s_2)=(-m,-n)}{\Res} \zq(s_1,s_2)=
\underset{(s_1,s_2)=(-m,-n)}{\Res}\zeta(s_1,s_2).
\end{equation}
\end{cor}
\begin{proof} By Thm.~\ref{thm:kn+1} and Table~\ref{Ta:dbzeta}
we only need to prove
\begin{multline}\label{equ:klenkey}
\lim_\qup (1-q)^{1-k}\sum_{r=0}^k (-1)^r {n+1-r\choose
k-r}{n\choose r}/(q^{n+1-r}-1) \\
=\begin{cases}
\displaystyle{\frac{-1}{n+1} } &\quad\text{if }k=0,\\
\displaystyle{(-1)^{1-k}\frac{B_k}{k}{n\choose k-1}}
&\quad\text{if }1\le k\le n.\end{cases}
\end{multline}
First by generating function of the Bernoulli numbers
\begin{equation*}
\frac{1}{q^{n+1-r}-1}=\frac{1}{e^{(n+1-r)\log q}-1}=
\sum_{l=0}^\infty \frac{B_l}{l!} \big((n+1-r)\log q\big)^{l-1}.
\end{equation*}
Plugging this into the left hand side of equation
\eqref{equ:klenkey}, replacing $1-q$ by $-\log q$, and exchanging
the summation we get
\begin{multline}\label{equ:exch}
(-\log q)^{1-k}\sum_{r=0}^k (-1)^r {n+1-r\choose k-r}{n\choose
r}/(q^{n+1-r}-1)\\
=(-1)^{1-k} \sum_{l=0}^\infty
\frac{B_l}{l!}(\log q)^{l-k} \sum_{r=0}^k (-1)^r {n+1-r\choose
k-r}{n\choose r} (n+1-r)^{l-1}.
\end{multline}
Then the inner sum over $r$ is the coefficient of $x^k$ of the
following polynomial
\begin{align}
f_l(x)=&\sum_{k=0}^{n+1} \sum_{r=0}^k (-1)^r {n+1-r\choose
k-r}{n\choose r} (n+1-r)^{l-1} x^k \nonumber \\
=&\sum_{r=0}^n (-1)^r {n\choose r} (n+1-r)^{l-1}
 \sum_{k=r}^{n+1} {n+1-r\choose k-r}x^k \nonumber\\
=&(x+1)^{n+1} \sum_{r=0}^n (-y)^r {n\choose r} (n+1-r)^{l-1},
\label{equ:auxf}
\end{align}
where $y=x/(x+1)$. When $l=0$ this expression becomes
\begin{equation*}
f_0(x)=\frac{(x+1)^{n+1}}{n+1}  \sum_{r=0}^n (-y)^r  {n+1 \choose
r} =\frac{(x+1)^{n+1}}{n+1}\Bigl[(1-y)^{n+1}-(-y)^{n+1}\Bigr]
=\frac{1}{n+1}\Bigl[1-(-x)^{n+1}\Bigr].
\end{equation*}
Note that $k\le n$ we see the coefficient of $x^k$ in $f_0(x)$ is
$0$ if $k>0$ and it's $1/(n+1)$ if $k=0$. If $k=0$ then only the
constant term $-1/(n+1)$ in \eqref{equ:exch} remains when $\qup$
which proves the corollary in this case. So we can assume $l,k>0$.
Then
\begin{align*}
f_l(x)=&\left(z\frac{d}{dz}\right)^{l-1}\left.\Bigl\{
(x+1)^{n+1}\sum_{r=0}^n (-y)^r
{n\choose r}z^{n+1-r}\Bigr\}\right|_{z=1}\\
=&\left(z\frac{d}{dz}\right)^{l-1}\left.\Bigl\{(x+1)^{n+1}z(z-y)^n
\Bigr\}\right|_{z=1}.
\end{align*}
Note that highest degree term in $f_l(x)$ is contained in
\begin{align*}
(x+1)^{n+1} \left.\left(\frac{d}{dz}\right)^{l-1}(z-y)^n
\right|_{z=1} =&n(n-1)\cdots (n-l+2) (x+1)^{n+1} (1-y)^{n-l+1}\\
=&n(n-1)\cdots (n-l+2) (x+1)^l.
\end{align*}
If $l=1$ one can easily modify this to get just $x+1$. If $l<k$
then the coefficient of $x^k$ in $f_l(x)$ is 0. If $l=k$ it is
equal to
$$n(n-1)\cdots (n-k+2)=(k-1)!{n\choose k-1}.$$
The last express is valid even for $k=l=1$. Thus the range of $l$
in the outer sum of \eqref{equ:exch} starts from $k$. Moreover,
the first term of \eqref{equ:exch} is
\begin{equation*}
(-1)^{1-k}\frac{B_k}{k}{n\choose k-1}
\end{equation*}
as desired. This completes the proof of the corollary.
\end{proof}

\begin{prop}
Let $k,n$ be two non-negative integers such that $n\ge k$ and $k$
is even. Let $m=k-n-2$. Then
\begin{equation*}
 \underset{(s_1,s_2)=(-m,-n)}{\Res} \zq(s_1,s_2)=
 \frac{-f(q)(q-1)/\log q} {D(q)},\quad D(q)=\prod_{j=n+1-k}^{n+1} F(q,j)^{\epsilon_j}
\end{equation*}
where $f(q)\in \Z[q]$ is a palindrome with leading coefficient
${n\choose k}$, $F(q,j)\in \Z[q]$ is a factor of $(q^j-1)/(q-1)$,
$\epsilon=0$ or $1$, and $\deg_q D(q)=n+\deg_q f(q)$, such
that
\begin{equation*}
\lim_\qup \underset{(s_1,s_2)=(-m,-n)}{\Res}
\zq(s_1,s_2)=\underset{(s_1,s_2)=(-m,-n)}{\Res} \zeta(s_1,s_2).
\end{equation*}
\end{prop}
\begin{proof} The computational proof is left as an
exercise for the interested readers.
\end{proof}

\begin{eg} By Thm.~\ref{thm:kn+1} we find with the help of Maple
\begin{equation*}
\underset{(s_1,s_2)=(4,-4)}{\Res}
\zq(s_1,s_2)=\frac{-2q^3(3q^2+4q+3)(q-1)/\log
q}{P_1(q,2)P_1(q,3)P_1(q,4)}.
\end{equation*}
where $P_a(q,m)=\sum_{j=0}^m  q^{aj}$. Moreover we can check that
\begin{equation*}
\lim_\qup \underset{(s_1,s_2)=(4,-4)}{\Res}
\zq(s_1,s_2)=\underset{(s_1,s_2)=(4,-4)}{\Res}
\zeta(s_1,s_2)=-\frac{1}{3}
\end{equation*}
by Table~\ref{Ta:dbzeta} with $k=2$ and $n=4$.
\end{eg}
\begin{eg} From Thm.~\ref{thm:kn+1} we get
\begin{equation*}
 \underset{(s_1,s_2)=(6,-8)}{\Res} \zq(s_1,s_2)=
 \frac{-14q^5g(q)(q-1)/\log q} {P_1(q,4)P_1(q,5)P_1(q,6)P_2(q,3)P_3(q,2)}
\end{equation*}
where $g(q)$ is a polynomial in $q$ of degree $14$ satisfying
$$q^{14}f(1/q)=f(q)=5q^{14}+6q^{13}+8q^{12}+7q^{11}-q^{10}-20q^9-30q^8-34q^7-\cdots.$$
Then we can compute with Maple
\begin{equation*}
\lim_\qup \underset{(s_1,s_2)=(6,-8)}{\Res}
\zq(s_1,s_2)=\underset{(s_1,s_2)=(6,-8)}{\Res}
\zeta(s_1,s_2)=\frac{7}{15}
\end{equation*}
by Table~\ref{Ta:dbzeta} with $k=4$ and $n=8$.
\end{eg}
\begin{eg} Consider the point $(s_1,s_2)=(5,-9)$. We have
\begin{equation*}
 \underset{(s_1,s_2)=(5,-9)}{\Res} \zq(s_1,s_2)=
 \frac{-42q^4g(q)(q-1)/\log q} {P_1(q,4)P_1(q,6)P_1(q,7)P_2(q,2)P_3(q,2)A(q,4)}
\end{equation*}
where $A(q,m)=\sum_{j=0}^m (-1)^j q^j$ and $g(q)$ is a polynomial
in $q$ of degree $18$ satisfying
$$q^{18}g(1/q)=g(q)=
2q^{18}-q^{17}-7q^{15}-11q^{14}-16q^{13}-4q^{12}+9q^{11}+28q^{10}+30q^9+\cdots
$$
so that we again have the equality
\begin{equation*}
\lim_\qup \underset{(s_1,s_2)=(5,-9)}{\Res}
\zq(s_1,s_2)=\underset{(s_1,s_2)=(5,-9)}{\Res}
\zeta(s_1,s_2)=-\frac{1}{2}
\end{equation*}
by Table~\ref{Ta:dbzeta} with $k=6$ and $n=9$.
\end{eg}

\begin{thm}\label{thm:kgen2}
Let $m,n$ be two non-negative integers and $k=m+n+2$. Then
$\zq(s_1,s_2)$ has indeterminacy at $(-m,-n)$ such that
\begin{align*}
\zq(-m,-n)=&(1-q)^{2-k}
\left\{\frac{(-1)^k}{(m+1)(n+1)(\log q)^2} \right.\\
+&\sum_{r=0}^{m}\frac{(-1)^{r+n+1}}{(n+1)\log q} {m\choose r}
\frac{1}{q^{m+1-r}-1}\\
+&\sum_{r=0}^n \frac{(-1)^{r+m+1}}{\log
q}\frac{m!(n+1-r)!}{(k-r)!} {n\choose
r}\frac{1}{q^{n+1-r}-1}\\
+& \left.\sum_{r_1=0}^{m}\sum_{r_2=0}^n (-1)^{r_1+r_2} {m\choose
r_1}{n\choose r_2} \frac{1}{q^{n+1-r_2}-1}
\frac{1}{q^{k-r_1-r_2}-1} \right\},
\end{align*}
and
\begin{align*}
\zq^R(-m,-n)=&(1-q)^{2-k} \left\{\sum_{r=0}^{m}
\frac{(-1)^{r+n+1}}{(n+1)\log q}
{m\choose r} \frac{1}{q^{m+1-r}-1}\right.\\
+& \sum_{r=0}^{m} \frac{(-1)^{r+n}}{\log q}{k-n-2\choose
r}\frac{n!(m+1-r)!}{(k-r)!}
\frac{q^{m+1-r}}{q^{m+1-r}-1}\\
+& \left.\sum_{r_1=0}^{m}\sum_{r_2=0}^n (-1)^{r_1+r_2} {m\choose
r_1}{n\choose r_2} \frac{1}{q^{n+1-r_2}-1}
\frac{1}{q^{k-r_1-r_2}-1}\right\}.
\end{align*}

\end{thm}
\begin{proof}
Use Thm.~\ref{thm:qmz}.
\end{proof}

Similar to Cor.~\ref{cor:compare1} and Cor.~\ref{cor:compare2} we have
\begin{cor} \label{cor:compare3}
Let $m$ and $n$ be two non-negative integers. Then
\begin{equation}\label{equ:zkgeqn}
\lim_\qup  \zq( -m,-n)= \zeta(s_1,s_2),\quad \lim_\qup  \zq^R(
-m,-n)= \zeta^R(s_1,s_2).
\end{equation}
\end{cor}
\begin{proof} Set $k=m+n+2$. We consider $\zq^R(-m,-n)$ first. From
Thm.~\ref{thm:kgen2} we get
\begin{equation*}
\zq^R(-m,-n)=(q-1)^{-k}\left(\frac{q-1}{\log q}\right)^2(A+B+C),
\end{equation*}
where
\begin{align*}
A=& \sum_{j=0}^\infty \frac{B_j}{j!} (\log q)^{j} \sum_{r=0}^{m}
\frac{(-1)^{m+r+1}}{n+1}
{m\choose r} (m+1-r)^{j-1}, \\
B=&\sum_{i=0}^\infty \frac{B_i}{i!} (\log q)^{i}
\sum_{r=0}^{m}(-1)^{r+m+i} {m\choose
r}\frac{n!(m+1-r)!}{(k-r)!} (m+1-r)^{i-1}, \\
C=&\sum_{i,j=0}^\infty \frac{B_iB_j}{i!j!}  (\log
q)^{i+j}\sum_{r_1=0}^{m}\sum_{r_2=0}^n (-1)^{k+r_1+r_2} {m\choose
r_1}{n\choose r_2} (n+1-r_2)^{i-1} (k-r_1-r_2)^{j-1}.
\end{align*}
We first compute $B$ as follows. Write
\begin{equation*}
B=\sum_{i=0}^\infty \frac{B_i}{i!} (\log q)^{i} W(m,n,i),
\end{equation*}
where
\begin{equation*}
W(m,n,i)=\sum_{r=0}^{m}(-1)^{r+k+n+i} \frac{m!n!}{r!(k-r)!}
(m+1-r)^i.
\end{equation*}
If $i=0$ then we can prove by decreasing induction on $n$ that
\begin{equation}\label{w0}
W(m,n,0)=\sum_{r=0}^{m}(-1)^{r+k+n}
\frac{m!n!}{r!(k-r)!}=\frac{1}{k(n+1)}.
\end{equation}
This is trivial if $n=k-2$. Suppose \eqref{w0} is true for $n\ge
1$ then we have
\begin{equation*}
W(m,n-1,0)=-\frac{k-n-1}{n}
W(m,n,0)+\frac{(n-1)!}{(n+1)!}=\frac{1}{kn},
\end{equation*}
as desired.

Similarly, we can compute $C$ as follows. Put
\begin{equation*}
C=\sum_{i,j=0}^\infty \frac{B_iB_j(\log q)^{i+j}}{i!j!(n+1)}
\sum_{r_1=0}^{m}\sum_{r_2=0}^n (-1)^{k+r_1+r_2} {m\choose
r_1}{n+1\choose r_2} (n+1-r_2)^{i} (k-r_1-r_2)^{j-1}.
\end{equation*}
We now change the upper limit of $r_2$ from $n$ to
$n+1$ in the above. The extra terms correspond to those
by setting $i=0$ and $r_2=n+1$, which produce exactly $A$. Therefore,
\begin{equation*}
C=\sum_{i,j=0}^\infty \frac{B_iB_j(\log
q)^{i+j}}{i!j!}V(m,n,i,j)-A
\end{equation*}
where
\begin{equation*}
V(m,n,i,j)=\sum_{r_1=0}^{m}\sum_{r_2=0}^{n+1}
\frac{(-1)^{k+r_1+r_2}}{n+1} {m\choose r_1}{n+1\choose r_2}
(n+1-r_2)^i (k-r_1-r_2)^{j-1}.
\end{equation*}
For $j\ge 1$ we have
\begin{align}
V(m,n,i,j)=&\left(x\frac{d}{dx}\right)^{j-1}\Bigg\{\left(y\frac{d}{dy}\right)^i
\Big\{\sum_{r_1=0}^{m}\sum_{r_2=0}^{n+1} \frac{(-1)^{k+r_1+r_2}}{n+1}\nonumber \\
\ &\phantom{some space just here} {m\choose r_1}{n+1\choose r_2}
y^{n+1-r_2}
x^{k-r_1-r_2}\Big\}\Big|_{y=1}\Bigg\}\Bigg|_{x=1}\nonumber\\
=&\frac{(-1)^k}{n+1}\left(x\frac{d}{dx}\right)^{j-1}\Bigg\{\left(y\frac{d}{dy}\right)^i
\Big\{ x(xy-1)^{n+1} (x-1)^{k-n-2}
 \Big\}\Big|_{y=1}\Bigg\}\Bigg|_{x=1}\label{equ:j=0}\\
=&\begin{cases} 0 &\text{ if $i+j<k$},\\
\displaystyle{\frac{(-1)^k n!}{(n+1-i)!}(k-i-1)!} \quad &\text{ if
}i+j=k, i\le n+1.\nonumber
\end{cases}
\end{align}
We have used the fact that if $i+j=k$ and $i> n+1$ then $j<k-n-1$
and by exchanging the two operators $x(d/dx)$ and $y(d/dy)$ we can
easily show that \eqref{equ:j=0} is zero. So if $l<k$ the total
contribution to the coefficient of $(\log q)^l$ from $V(m,n,i,j)$
with $j>0$ is trivial and if $l=k$ it is equal to
\begin{equation}\label{equ:midterms}
\sum_{i=0}^{n+1} B_iB_{k-i}\frac{(-1)^k n!(k-i-1)!}{i!(k-i)!(n+1-i)!}
=\begin{cases}
\displaystyle{\frac{B_{k-1}}{2(k-1)}} &\text{if } 2\nmid k,\\
\displaystyle{\frac{B_k}{k(n+1)}+ \sum_{i=1}^{n+1} \frac{B_iB_{k-i}}{i(k-i)}
{n\choose i-1}} &\text{if } 2\mid k,
\end{cases}
\end{equation}
because $ k\ge n+2\ge 2$ and $B_k=0$ if $k$ is odd.

To deal with $V(m,n,i,0)$ note that equation~\eqref{equ:j=0} still
makes sense if we interpret the operator $x(d/dx)^{-1}$ as
follows:
\begin{equation*}
\left(x\frac{d}{dx}\right)^{-1}\Big\{F(x)\Big\}\Big|_{x=1}=\int_0^1\frac{F(x)}{x} dx
\end{equation*}
whenever $F(0)=0$. Thus we get
\begin{equation*}
V(m,n,i,0)= \frac{(-1)^{k}}{n+1}\left(y\frac{d}{dy}\right)^i
\Bigg\{ \int_0^1 (xy-1)^{n+1} (x-1)^{m}\,dx\Bigg\} \Bigg|_{y=1}.
\end{equation*}
Therefore if $i=0$ then we get
\begin{equation}\label{equ:00cancel}
V(m,n,0,0)=\frac{(-1)^{k}}{n+1} \int_0^1
(x-1)^{k-1}\,dx=-\frac{1}{k(n+1)}=-W(m,n,0)
\end{equation}
from equation~\eqref{w0}. If $i\ge 1$ then integrating by parts we
get
\begin{align*}
\ & \int_0^1 (xy-1)^{n+1} (x-1)^{m}\,dx\\
=&\frac{(xy-1)^{n+2}}{y(n+2)}(x-1)^{m}\Bigg|_0^1-\frac{m}{y(n+2)}\int_0^1
(xy-1)^{n+2}(x-1)^{m-1}\,dx\\
= & \cdots \cdots \\
=&(-1)^{k+1}\left(\frac{1}{y(n+2)}-\frac{m}{y^2(n+2)(n+3)}+\cdots
+(-1)^{m-1}\frac{m!(n+1)!}{y^{m}(m+n+1)!}\right. \\
\ &\ \qquad \left. +(-1)^{m}\frac{m!(n+1)!}{y^{m}(m+n+1)!}\int_0^1
(xy-1)^{m+n+1}\,dx, \right)\\
=&(-1)^{k+n}\frac{m!(n+1)!}{k!}
\frac{(y-1)^k}{y^{m+1}}+\sum_{r=0}^{m}(-1)^{r+k+1}
\frac{m!(n+1)!}{y^{r+1}(m-r)!(n+2+r)!} .
\end{align*}
It follows from changing the index $r$ to $m-r$ that
\begin{align}
V(m,n,i,0)=&(-1)^{k+n} \frac{m!n!}{k!}
\left(y\frac{d}{dy}\right)^i \Bigg\{
\frac{(y-1)^k}{y^{m+1}}\Bigg\}
\Bigg|_{y=1}+\sum_{r=0}^{m}(-1)^{r+n+i+1}
\frac{m!n!}{r!(k-r)!}(m+1-r)^i \nonumber\\
 =&\begin{cases}
(-1)^{k+1} W(m,n,0) \ &\text{ if } 0<i<k,\\
(-1)^{k+1} W(m,n,0)+(-1)^{k+n} n!(k-n-2)! \quad &\text{ if
}i=k.\label{equ:dangle}
\end{cases}
\end{align}
Thus when $0<i<k$ and $k$ is even we have $V(m,n,i,0)=-W(m,n,i)$.
It follows from \eqref{equ:midterms}, \eqref{equ:00cancel} and
\eqref{equ:dangle} that
\begin{equation*}
\lim_\qup\zq^R(-m,-n)= \zeta^R(s_1,s_2)
\end{equation*}
since $B_k=0$ if $k>2$ is odd.

Let's turn to prove the first equality in \eqref{equ:zkgeqn}. By
Thm.~\ref{thm:kgen2} we have
\begin{equation*}
\zq(-m,-n)=(q-1)^{-k}\left(\frac{q-1}{\log q}\right)^2(D+A+E+C)
\end{equation*}
where $A$ and $C$ are as above and
\begin{equation*}
D=\frac{1}{(m+1)(n+1)},\quad E=\sum_{i=0}^\infty \frac{B_i}{i!}
(\log q)^{i} U(m,n,i),
\end{equation*}
where
$$U(m,n,i)=-\sum_{r=0}^n(-1)^{r+n}\frac{m!n!}{r!(k-r)!}(n+1-r)^i.$$
Hence
\begin{equation}\label{equ:u0}
U(m,n,0)=W(n,m,0)=\frac{-1}{k(m+1)}=\frac{-1}{k(k-n-1)}=\frac{1}{k(n+1)}-D.
\end{equation}
We only need to show that
\begin{equation}\label{equ:u=w}
U(m,n,i)=
\begin{cases} W(m,n,i) &\text{ if }0<i<k,\\
W(m,n,i)-(-1)^{n} n!(k-n-2)! \quad &\text{ if }i=k.
\end{cases}
\end{equation}
Indeed when $i>0$ we have
\begin{align*}
U(m,n,i)=& (-1)^{n+1}\left(y\frac{d}{dy}\right)^i \Bigg\{
\sum_{r=0}^n (-1)^r \frac{m!n!}{k!}{k\choose r} y^{n+1-r} \Bigg\}\Bigg|_{y=1} \\
=& (-1)^{n+1}\frac{m!n!}{k!}\left(y\frac{d}{dy}\right)^i
\Bigg\{y^{n+1-k}(y-1)^k
-\sum_{r=n+1}^k (-1)^r {k\choose r} y^{n+1-r}\Bigg\}\Bigg|_{y=1} \\
=& (-1)^{n+1}\frac{m!n!}{k!} \left(y\frac{d}{dy}\right)^i \Bigg\{
\frac{(y-1)^k}{y^{m+1}}\Bigg\} +\sum_{r=n+1}^k
(-1)^{r+n}\frac{m!n!}{r!(k-r)!}(n+1-r)^i.
\end{align*}
Then first term is 0 if $i<k$ and it's $(-1)^{n+1}m!n!$ if $i=k$.
When $r=n+1$ the summand in the second term is zero since $i>0$.
So we can let $r$ range only from $n+2$ to $k$. Then change the
index $r$ to $k-r$ (and let $r$ run from $0$ to $m$) we can see
immediately that the second term is exactly $W(m,n,i)$. This
proves equation \eqref{equ:u=w} which together with \eqref{equ:u0}
implies the first equation in \eqref{equ:zkgeqn}. We thus finish
the proof of the corollary.
\end{proof}

We conclude this section by remarking that by shuffle relation
\eqref{equ:dbshf} we can also analyze $\zq(s_1,s_2)$ at
$(-n,n+2-k)$ for any non-negative integers $k$ and $n$. For
example, it's easy to compute directly that
\begin{align*}
\underset{(s_1,s_2)=(-3,2)}{\Res}
\zq(s_1,s_2)=&-\frac{1}{(1-q)\log q}\sum_{r=0}^3 (-1)^r{3\choose
r}(r+1)\frac{q^{r+1}}{1-q^{r+1}}\\
=&\frac{-q(q-1)^2}{(q+1)(q^2+1)(q^2+q+1)\log q},
\end{align*}
which can be obtained also by the shuffle relation \eqref{equ:dbshf}
and the expression
\begin{equation*}
\underset{(s_1,s_2)=(2,-3)}{\Res} \zq(s_1,s_2)
=\frac{q(q-1)^2}{(q+1)(q^2+1)(q^2+q+1)\log q}
\end{equation*}
by taking $k=n=3$ in Thm.~\ref{thm:kn+1}. Thus $(-3,2)$ is a
simple pole of $\zq(s_1,s_2)$. On the other hand $\zeta(s_1,s_2)$
does not have a pole along $s_1+s_2=-1$. Indeed we find that
$$\lim_\qup \underset{(s_1,s_2)=(-3,2)}{\Res}\big\{\zq(s_1,s_2)\big\}=0.$$

\section{$q$-multiple polylogarithms}\label{sec:polylog}
It is well known that special values of the multiple zeta function
$\zeta(s_1,\dots,s_d)$ at positive integers $(n_1,\dots,n_d)$ can
be regarded as single-valued version of multiple polylogarithm
$Li_{n_1,\dots,n_d}(z_1,\dots,z_d)$ evaluated at
$z_1=\dots=z_d=1$. For $|z_j|<1$ these functions can be defined as
\begin{equation*}
 Li_{n_1,\dots,n_d}(z_1,\dots,z_d)=
 \sum_{0<k_1<\dots<d_d}\frac{z_1^{k_1}\dots z_d^{k_d}}
{k_1^{n_1}\dots k_d^{n_d}}.
\end{equation*}
By Chen's iterated integral
\begin{equation}\label{chen}
(-1)^d\int_0^1
\frac{dt_1}{t_1-a_1}\circ\left(\frac{dt_1}{t_1}\right)^{\circ
(n_1-1)} \circ\dots \circ
\frac{dt_d}{t_d-a_d}\circ\left(\frac{dt_d}{t_d}\right)^{\circ
(n_d-1)},
\end{equation}
where $a_j=1/\prod_{i=j}^d z_i$ for all $j=1,\dots,d$, we can
obtain the analytic continuation of this function as
a multi-valued function on $\C^d\setminus\DD_d$ where
\begin{equation*}
\DD_d=\Bigl\{(z_1,\dots,z_d)\in \C^d:  \prod_{i=j}^d z_i=1,\
j=1,\dots,d\Bigr\}.
\end{equation*}
When $|z_j|<$ we define its $q$-analog ($0<q<1$) by
\begin{equation*}
 Li_{q;n_1,\dots,n_d}(z_1,\dots,z_d)=
 \sum_{0<k_1<\dots<d_d}\frac{z_1^{k_1}\dots z_d^{k_d}}
{[k_1]^{n_1}\dots [k_d]^{n_d}}.
\end{equation*}
Clearly when $\qup$ we recover the ordinary multiple
polylogarithm. Moreover, the special value of $q$-multiple zeta
function at positive integers
\begin{equation*}
\zq(n_1,\dots,n_d)=
Li_{q;n_1,\dots,n_d}\big(q^{n_1-1},\dots,q^{n_d-1}\big).
\end{equation*}

Note that our definition of the $q$-multiple polylogarithms is
different from that of \cite{Sch}. In case of logarithm and dilogarithm
our definitions are different from that of \cite{K}.
We want to convince the readers
that ours are also good analogs of the ordinary ones.

We can mimic the method in section~\ref{sec:fqzq} to
get the analytic continuation of $Li_{q;n_1,\dots,n_d}(z_1,\dots,z_d)$.
\begin{thm}\label{thm:anaqlog}
The $q$-multiple polylogarithm function
$Li_{q;n_1,\dots,n_d}(z_1,\dots,z_d)$ converges if $|z_j|<1$ for
all $j=1,\dots,d$. It can be analytically continued to a
multi-valued function over $\C^d\setminus \DD_{q;d}$ via the
series expansion
\begin{equation}\label{qlogmerom}
Li_{q;n_1,\dots,n_d}(z_1,\dots,z_d)=(1-q)^{n_1+\dots+n_d}
\sum_{r_1,\dots,r_d=0}^{+\infty} \prod_{j=1}^d
\left[{n_j+r_j-1\choose r_j} \frac{z_j^jq^{jr_j}}{1-(z_j \cdots
z_d)q^{r_j+\cdots+r_d}} \right].
\end{equation}
\end{thm}
\begin{proof}
The first part of the lemma is obvious.
Let's concentrate on the analytic continuation. By binomial expansion
$(1-x)^{-n}=\sum_{r=0}^{\infty}{n+r-1\choose r} x^r$ we get
 \begin{equation*}
Li_{q;n_1,\dots,n_d}(z_1,\dots,z_d)=(1-q)^{n_1+\dots+n_d}
\sum_{0<k_1<\dots<k_d} \sum_{r_1,\dots,r_d=0}^{+\infty}
 \prod_{j=1}^d {n_j+r_j-1\choose r_j} (z_jq^{r_j})^{k_j}.
\end{equation*}
As $0<q<1$ the series converges absolutely by Stirling's formula
so we can exchange the summations. The theorem follows immediately
from Lemma~\ref{lem:uuu} by taking $x_j=z_jq^{r_j}$.
\end{proof}

However, this analytic continuation is not suitable for comparing
with its ordinary counterpart. The we define Jackson's
$q$-differential operator (cf. \cite{J1}) by
\begin{equation*}
D_{q;z} f(z)= \frac{f(z)-f(qz)}{(1-q)z}.
\end{equation*}
\begin{lem} \label{lem:qdiff}
Let $d,n_1,\dots,n_d$ be positive integers.
If $n_j\ge 2$ then we have
\begin{equation*}
    D_{q;z_j}Li_{q;n_1,\dots,n_d}(z_1,\dots,z_d)=
\frac{1}{z_j}Li_{q;n_1,\dots,n_j-1,\dots,n_d}(z_1,\dots,z_d);
\end{equation*}
if $d\ge 2$ and $n_j=1$ then
\begin{multline*}
D_{q;z_j}Li_{q;n_1,\dots,n_d}(z_1,\dots,z_d)
=\frac{1}{1-z_j}Li_{q;n_1,\dots,\widehat{n_j},\dots,n_d}(z_1,\dots,z_{j-1}z_j,\dots,z_d)
    \\
-\frac{1}{z_j(1-z_j)}Li_{q;n_1,\dots,\widehat{n_j},\dots,n_d}(z_1,\dots,z_jz_{j+1},\dots,z_d).
\end{multline*}
Here the second term does not appear if $j=d$. If $d=n_1=1$ then
\begin{equation*}
    D_{q;z}Li_{q;1}(z)=\frac{1}{1-z}.
\end{equation*}
\end{lem}
\begin{proof} Clear.
\end{proof}
The same properties listed in the lemma are satisfied by the
ordinary multiple polylogarithms. We note that the first equation
in \cite[Lemma 1]{Sch} is valid only for $n_j\ge 2$.

Recall that for any continuous function $f(x)$ on $[a,b]$
Jackson's $q$-integral (cf. \cite{J2}) is defined by
\begin{equation*}
\int_a^b f(x) \,d_qx:=\sum_{i=0}^\infty f\big(a+q^i(b-a)\big)
(q^i-q^{i+1})(b-a).
\end{equation*}
Then for every $x\in [a,b]$ we have
\begin{equation}\label{qintconv}
\int_a^x D_{q;t} f(t) \,d_qt= f(x)-f(a),\quad \text{and }\quad
\lim_\qup \int_a^x  f(t) \,d_qt=  \int_a^x  f(t) \,dt.
\end{equation}

\begin{rem} Note that in general $\int_a^b f(x) \,d_qx+\int_b^c f(x)
\,d_qx\ne \int_a^c f(x) \,d_qx.$
\end{rem}

Similar to Chen's iterated integrals we can define the
$q$-iterated integrals as follows:
\begin{equation*}
\int_a^b \frac{d_qt_1}{t_1-a_1}\circ \dots \circ
\frac{d_qt_r}{t_r-a_r}:=\int_a^b\left(\int_a^{t_r} \cdots
\int_a^{t_3} \left(\int_a^{t_2}
\frac{d_qt_1}{t_1-a_1}\right)\frac{d_qt_2}{t_2-a_2}  \dots
\frac{d_qt_{r-1}}{t_{r-1}-a_{r-1}}\right)\frac{d_qt_r}{t_r-a_r}.
\end{equation*}
Define
\begin{equation*}
\DD_{q;d}=:\Bigl\{(z_1,\dots,z_d)\in \C^d: \prod_{i=j}^d
z_i=q^{-m},\ m\in \Z_{\ge0},\ j=1,\dots,d \Bigr\}.
\end{equation*}

\begin{cor}
We can analytically continue $Li_{q;n_1,\dots,n_d}(z_1,\dots,z_d)$
to $\C^d\setminus \DD_{q;d}$ by the $q$-iterated integral
\begin{equation*}
(-1)^d\int_0^1
\frac{d_qt_1}{t_1-a_1}\circ\left(\frac{d_qt_1}{t_1}\right)^{\circ
(n_1-1)} \circ\dots \circ
\frac{d_qt_d}{t_d-a_d}\circ\left(\frac{d_qt_d}{t_d}\right)^{\circ
(n_d-1)},
\end{equation*}
where $a_j=1/\prod_{i=j}^d z_i$ for all $j=1,\dots,d$.
Further, for all $(z_1,\dots,z_d)\in \C^d$ such that
$\prod_{i=j}^d z_i\notin [1,+\infty)$ for $1\le j\le d$ we have
$$\lim_\qup Li_{q;n_1,\dots,n_d}(z_1,\dots,z_d)
=Li_{n_1,\dots,n_d}(z_1,\dots,z_d),$$
where $Li_{n_1,\dots,n_d}(z_1,\dots,z_d)$ is defined by \eqref{chen}
with the path being the straight line segment from $0$ to $1$ in $\C^1$.
\end{cor}
\begin{proof}
It follows from Lemma~\ref{lem:qdiff} and equation
\eqref{qintconv}. The singular set $\DD_{q;d}$ is determined by
Thm.~\ref{thm:anaqlog} so that for each $j=1,\dots,d$ the function
$1/(t-a_j)$ is continuous on $[0,1]$.
\end{proof}

\section{Iterated integral $q$-shuffle relations}\label{sec:shuffle2}
In section \ref{sec:shuffle1} we encountered some $q$-shuffle
relations of $q$-multiple zeta functions. Classically, multiple
zeta  values satisfy another kind of shuffle relation coming from
their representations by Chen's iterated integrals. In our setting
we have seen that special values of $q$-multiple zeta functions
can be also represented by $q$-iterated integrals. In this last
section we would like to study the shuffle relations related to
these $q$-iterated integrals. We shall see that they're more
involved than their ordinary counterparts. We start by writing
\begin{equation*}
\sh(u_1\circ \cdots \circ  u_r , u_{r+1}\circ \cdots \circ u_{r+s}
)=\sum_{\gs } u_{\gs(1)}\circ \cdots   \circ u_{\gs(r+s)}.
\end{equation*}
where $\gs$ runs through all the permutations of $\{1,\dots,r+s\}$
such that $\gs^{-1}(a)<\gs^{-1}(b)$ whenever $1\le a<b\le r$ or
$r+1\le a<b\le r+s$. For any expressions $F_i$ we put
$$\bigsqcup_{i=1}^r F_i= F_1\circ \cdots \circ  F_r$$

\begin{lem}\label{lem:qshfl}
Let $u_i=d_qt/(t-a_i)$ and $v_j=d_qt/(t-b_i)$ where
$|a_i|,|b_j|\le 1$ for all $1\le i\le r$ and $1\le j\le s$. Let
$a$ be any positive number. Then
\begin{multline}
\int_0^a u_1\circ \cdots \circ  u_r  \cdot \int_0^a v_1\circ
\cdots \circ  v_s=  \int_0^a  \sh(u_1\circ \cdots \circ  u_r
, v_1\circ \cdots \circ  v_s )+ \sum_{c=1}^{\min(r,s)} (q-1)^c\cdot \\
 \cdot \sum_{\substack{1\le i_1<\dots<i_c\le r\\
1\le j_1<\dots<j_c\le s}}\int_0^a \bigsqcup_{\ga=1}^{c+1}
\sh(u_{1+i_{\ga-1}}\circ \cdots \circ u_{i_\ga-1},
v_{1+j_{\ga-1}}\circ \cdots \circ v_{j_\ga-1}) \circ\langle
u_{i_\ga}, v_{j_\ga}\rangle, \label{equ:shf}
\end{multline}
where $i_0=j_0=0$, $i_{c+1}=r+1$, $j_{c+1}=s+1$, $\langle u_{r+1},
v_{s+1}\rangle=1$, and for all $i,j$
\begin{equation*}
\langle u_i, v_j\rangle = \frac{t d_qt}{(t-a_i)(t-b_j)}=
\begin{cases}
\frac{1}{b_j-a_i}\left(\frac{b_jd_qt}{t-b_j}-\frac{a_id_qt}{t-a_i}\right)
\quad&\text{ if }a_i\ne b_j,\\
\frac{d_qt}{t-b}+\frac{b d_qt}{(t-b)^2} &\text{ if }a_i=b_j=b.
\end{cases}
\end{equation*}
\end{lem}
\begin{proof} This can be proved by induction on $r+s$.
We only want to mention that the key formula is
\begin{equation*}
D_{q;x} [f(x)g(x)]=  [D_{q;x}f(x)]g(x)+ f(x) [D_{q;x} g(x)]+
x(q-1) [D_{q;x} f(x)][D_{q;x} g(x)].
\end{equation*}
\end{proof}

We will say the term $\langle u_i, v_j\rangle$ is a {\em collapse}
in the shuffle. The lemma roughly says that $q$-iterated shuffle
relations is different from those produced by Chen's iterated
integrals because collapses may occur. The number of collapses is
at most $\min(r,s)$.

Lemma \ref{lem:qshfl} implies that if $m, n\ge 2$ are different
then
\begin{align}\nonumber
\zq(m)\zq(n)=&\int_0^1\frac{d_qt}{t-q^{1-m}}\circ
\left(\frac{d_qt}{t}\right)^{\circ (m-1)}\cdot \int_0^1
\frac{d_qt}{t-q^{1-n}}\circ\left(\frac{d_qt}{t}\right)^{\circ(n-1)} \\
=&A_q(m,n)+A_q(n,m)+B_q(m,n),\label{ABC}
\end{align}
where
\begin{align*}
A_q(m,n)=&\sum_{a=0}^{m-1}\sum_{c=0}^{\min(a,n)}
E(a,n;c)\int_0^1\frac{d_qt}{t-q^{1-m}}\circ
\left(\frac{d_qt}{t}\right)^{\circ
(m-1-a)}\circ\frac{d_qt}{t-q^{1-n}}
\circ\left(\frac{d_qt}{t}\right)^{\circ(n+a-1-c)},\\
B_q(m,n)=&(q-1)\sum_{c=0}^{\min(m,n)-1} E(m-1,n-1;c)
\int_0^1\frac{td_qt}{(t-q^{1-m})(t-q^{1-n})}\circ
\left(\frac{d_qt}{t}\right)^{\circ (m+n-2-c)}.
\end{align*}
Here the coefficient $E(r,s;c)$ represents $(q-1)^c$ times the
numbers of terms in shuffle of $u_1\circ \cdots \circ u_r$ and
$v_1\circ \cdots \circ v_s$ with $c$ collapses (see
\eqref{equ:shf}). It is not hard to see that $E(r,s;0)={r+s\choose
r}$. Thus
\begin{equation*}
E(r,s;c)=(q-1)^c \sum_{\substack{1\le i_1<\dots<i_c\le r\\
1\le j_1<\dots<j_c\le s}} \prod_{\ga=1}^{c+1}
{i_\ga+j_\ga-i_{\ga-1}-j_{\ga-1}-2 \choose i_\ga-i_{\ga-1}-1}.
\end{equation*}
We want to convert the expressions in \eqref{ABC} into something
that is close to linear combinations of multiple zeta functions.
Since $m\ne n$ we get
\begin{align*}
\ &\int_0^1\frac{td_qt}{(t-q^{1-m})(t-q^{1-n})}\circ
\left(\frac{d_qt}{t}\right)^{\circ (m+n-2-c)}\\
=&\int_0^1\left(\frac{1}{1-q^{m-n}}\frac{d_qt}{t-q^{1-m}}
+\frac{1}{1-q^{n-m}}\frac{d_qt}{t-q^{1-n}}\right)\circ
\left(\frac{d_qt}{t}\right)^{\circ (m+n-2-c)}\\
=&\frac{1}{q^{m-n}-1}\sum_{k=1}^\infty
\frac{q^{(m-1)k}}{[k]^{m+n-1-c}}+\frac{1}{q^{n-m}-1}\sum_{k=1}^\infty
\frac{q^{(n-1)k}}{[k]^{m+n-1-c}}.
\end{align*}
\begin{prop} \label{prop:Bqmn}
For any positive integers $m,n$ we have
\begin{multline*}
B_q(m,n)=(q-1)\sum_{c=0}^{\min(m,n)-1}
E(m-1,n-1;c) \cdot  \\
\left(\frac{1}{q^{m-n}-1}\sif^{n-1-c}\zq(m+n-1-c)
+\frac{1}{q^{n-m}-1}\sif^{m-1-c}\zq(m+n-1-c)\right).
\end{multline*}
\end{prop}
\begin{proof} It follows from Cor.~\ref{cor:shift}.
\end{proof}
To handle $A_q(m,n)$ we need to evaluate
\begin{equation*}
Li_{q;m-\ga,\gb}(q^{m-n},q^{n-1})=\int_0^1\frac{d_qt}{t-q^{1-m}}\circ
\left(\frac{d_qt}{t}\right)^{\circ
(m-\ga-1)}\circ\frac{d_qt}{t-q^{1-n}}
\circ\left(\frac{d_qt}{t}\right)^{\circ(\gb-1)}
\end{equation*}
where $0\le \ga\le m-1$ and $\max(n,\ga)\le \gb\le n+\ga$. By
Cor.~\ref{cor:shift} we get
\begin{equation}\label{ligb}
Li_{q;m-\ga,\gb}(q^{m-n},q^{n-1})=\sum_{i=0}^{\gb-n} {\gb-n\choose
i}(1-q)^i Li_{q;m-\ga,\gb-i}(q^{m-n},q^{\gb-i-1}).
\end{equation}
So we need to evaluate
\begin{equation*}
Li_{q;m-\ga,\gam}(q^{m-n},q^{\gam-1})=\sum_{1\le k<l}
\frac{q^{(m-n)k}}{[k]^{m-\ga}}\frac{q^{(\gam-1)l}}{[l]^\gam}.
\end{equation*}
If $n>\ga$ then by Cor.~\ref{cor:shift} we get
\begin{equation}\label{ligam}
Li_{q;m-\ga,\gam}(q^{m-n},q^{\gam-1})=\sum_{j=0}^{n-\ga-1}
{n-\ga-1\choose j}(1-q)^j \zq(m-\ga-j, \gam).
\end{equation}
By taking $\ga=a$, $\gb=n+a-c$ and $\gam=\gb-i$ in \eqref{ligb}
and \eqref{ligam} we get
\begin{prop}\label{prop:n>m} If $n>m$ then
\begin{equation*}
A_q(m,n)=\sum_{a=0}^{m-1}\sum_{c=0}^{\min(a,n)}\sum_{i=0}^{a-c}\sum_{j=0}^{n-a-1}
E(a,n;c) {a-c\choose i}(1-q)^{i+j} {n-a-1\choose j} \zq(m-a-j,
n+a-c-i).
\end{equation*}
\end{prop}

We now consider the case $n<m$. For $j>0$ define
\begin{equation*}
\caT^j\zq(\gam)=\sum_{l\ge 1}
\frac{q^j-q^{jl}}{1-q^j}\frac{q^{(\gam-1)l}}{[l]^\gam}.
\end{equation*}
Set $\caT^0\zq(\gam)=\lim_{j\to 0} \caT^j\zq(\gam)=\sum_{l\ge 1}
(l-1) q^{(\gam-1)l}/[l]^\gam.$ To evaluate these we need
\begin{lem} For any $e\ge 0$ we have
\begin{equation*}
Li_{q;\gam}(q^{e+\gam})=\sum_{l\ge 1}
\frac{q^{(e+\gam)l}}{[l]^\gam}=\sum_{i=0}^{\gam-2} (q-1)^i
{i+e\choose e} \zq(\gam-i)+(q-1)^{\gam-2} \sum_{l\ge 1}
\frac{q^l(q^{(e+1)l}-1)}{[l]^2}.
\end{equation*}
\end{lem}
\begin{proof} Set $a(e,\gam)=Li_{q;\gam}(q^{e+\gam})$. If $e=0$ we have
\begin{equation*}
a(0,\gam)=Li_{q;\gam}(q^{\gam})=\sum_{i=0}^{\gam-2} (q-1)^i
\zq(\gam-i)+
 (q-1)^{\gam-1} \sum_{l\ge 1} \frac{q^{l}}{[l]},
\end{equation*}
which confirms the lemma in this case. If $e\ge 1$ we have
\begin{align*}
a(e,\gam)=&\sum_{j=0}^{\gam-2} (q-1)^j a(e-1,\gam-j) +
(q-1)^{\gam-1} \sum_{l\ge 1} \frac{q^{(e+1)l}}{[l]}\\
=&\sum_{j=0}^{\gam-2} (q-1)^j \sum_{i=0}^{\gam-j-2} (q-1)^i
{i+e-1\choose e-1} \zq(\gam-i-j) + (q-1)^{\gam-1} \sum_{l\ge 1}
\sum_{i=0}^e  \frac{q^{(i+1)l}}{[l]}\\
=&\sum_{\gs=0}^{\gam-2} (q-1)^\gs {\gs+e\choose e}
\zq(\gam-\gs)+(q-1)^{\gam-2} \sum_{l\ge 1}
\frac{q^l(q^{(e+1)l}-1)}{[l]^2}
\end{align*}
by induction and the combinatorial identity
\begin{equation}\label{combid}
\sum_{i=0}^\gs {i+e-1\choose e-1}={\gs+e\choose e}.
\end{equation}
\end{proof}

\begin{cor} \label{cor:Tjzq}
For all $j\ge 1$ we have
\begin{equation*}
\caT^j\zq(\gam)=-\zq(\gam)+\sum_{i=1}^{\gam-2}
\frac{(q-1)^i}{q^j-1} {i+j-1\choose j-1} \zq(\gam-i)+
\frac{(q-1)^{\gam-2}}{q^j-1}\sum_{l\ge
1}\frac{q^l(q^{jl}-1)}{[l]^2}.
\end{equation*}
\end{cor}

Now we fix $\gam$ and for all $s\ge r\ge 1$ we set
\begin{equation*}
X_\gam(r,s):=X(r,s):=Li_{q;r,\gam}(q^s,q^{\gam-1}).
\end{equation*}
\begin{lem} \label{lem:X} For every $e\ge 0$ we have
\begin{equation*}
X_\gam(r,r+e)=\sum_{i=0}^{r-1} (q-1)^i {i+e\choose e}
\zq(r-i,\gam)+(q-1)^r \sum_{j=0}^e \caT^j\zq(\gam).
\end{equation*}
\end{lem}
\begin{proof}
Clearly
\begin{align*}
X(r,r)=&\zq(r,\gam)+\sum_{1\le k<l}
\frac{q^{(r-1)k}(q^k-1)}{[k]^r}\frac{q^{(\gam-1)l}}{[l]^\gam}\\
=&\zq(r,\gam)+(q-1)X(r-1,r-1)\\
=&\sum_{i=0}^{r-1} (q-1)^i \zq(r-i,\gam)+(q-1)^r\sum_{l\ge 1}
\frac{(l-1)q^{(\gam-1)l}}{[l]^\gam}.
\end{align*}
And for $e\ge 1$
\begin{align*}
X(r,r+e)=&X(r,r+e-1)+(q-1)X(r-1,r+e-1) \\
=&\sum_{i=0}^{r-1} (q-1)^i X(r-i,r+e-1-i)+(q-1)^r  \sum_{1\le k<l}
q^{ek}\frac{q^{(\gam-1)l}}{[l]^\gam}\\
=&\sum_{i=0}^{r-1} (q-1)^i X(r-i,r+e-1-i)+(q-1)^r  \sum_{l\ge 1}
\frac{q^e-q^{el}}{1-q^e}\frac{q^{(\gam-1)l}}{[l]^\gam}\\
=&\sum_{i=0}^{r-1} (q-1)^i \left(\sum_{j=0}^{r-i-1} {j+e-1\choose
e-1}(q-1)^j \zq(r-i-j,\gam) \right)+(q-1)^r \sum_{j=0}^e
\caT^j\zq(\gam)
\end{align*}
by induction. Let $\gs=i+j$ in the first sum of the last
expression. Then the lemma follows from \eqref{combid}.
\end{proof}

Define $\xi_q(j)=\sum_{l\ge 1}q^{(j+1)l}/[l]^2.$ We now have
\begin{prop} \label{prop:n<m} If $n<m$ then
\begin{align*}
A_q(m,n)=&\sum_{a=0}^{n-1}\sum_{c=0}^{\min(a,n)}\sum_{i=0}^{a-c}\sum_{j=0}^{n-a-1}
E(a,n;c) {a-c\choose i}(1-q)^{i+j} {n-a-1\choose j} \zq(m-a-j,
n+a-c-i)\\
+&\sum_{a=n}^{m-1}\sum_{c=0}^{\min(a,n)}\sum_{i=0}^{a-c}E(a,n;c){a-c\choose
i}(1-q)^{i} X_{n+a-c-i}(m-a,m-n),
\end{align*}
where
\begin{align*}
X_\gam(r,s)=&\sum_{i=0}^{r-1} (q-1)^i {i+s-r\choose s-r}
\zq(r-i,\gam)+(q-1)^r \sum_{j=0}^{s-r}  \sum_{i=1}^{\gam-2}
\frac{(q-1)^i}{q^j-1} {i+j-1\choose j-1} \zq(\gam-i) \\
+& (q-1)^r \sum_{j=0}^{s-r}
\frac{(q-1)^{\gam-2}}{q^j-1}\big(\xi_q(j)-\zq(2)\big)-(q-1)^r
(s-r+1)\zq(\gam) .
\end{align*}
\end{prop}
\begin{proof} This is clear from equations \eqref{ligb}
and \eqref{ligam}, Lemma~\ref{lem:X} and Cor.~\ref{cor:Tjzq}.
\end{proof}

Putting everything together we arrive at
\begin{thm}
Let $m\ne n$ be two positive integers no less than 2. Then
\begin{equation*}
    \zq(m)\zq(n)=A_q(m,n)+A_q(n,m)+B_q(m,n),
\end{equation*}
where $A_q(m,n)$ is given by Prop.~\ref{prop:n>m} and
Prop.~\ref{prop:n<m}, and $B_q(m,n)$ is given by
Prop.~\ref{prop:Bqmn}.
\end{thm}
It is not hard to see that when $\qup$ we recover the ordinary
shuffle relations of the MZVs originally produced by using Chen's
iterated integrals. The only unpleasant terms in our $q$-analogs
are given by $\xi_q(j)$ which is closely related to $\zq(2)$.

\noindent {\em Email:} jqz@math.upenn.edu

\end{document}